\documentclass[11pt,a4paper]{article}
\usepackage[a4paper]{geometry}
\usepackage{amssymb,latexsym,amsmath,amsfonts,amsthm}
\usepackage{graphicx}
\usepackage{algorithm,algorithmic}
\usepackage{dsfont}
\usepackage{epstopdf}
\usepackage{subfigure}
\usepackage{mathrsfs,enumerate}
\usepackage{tikz}
\usepackage{epsfig}
\usepackage{booktabs}
\usepackage{mathtools}
\usepackage{comment,verbatim}
\usepackage{overpic}
\usepackage{hyperref}
\usepackage{bbm}
\usepackage[toc,page]{appendix}

\newtheorem{theorem}{Theorem}[section]
\newtheorem{lemma}[theorem]{Lemma}

\theoremstyle{definition}

\theoremstyle{definition}

\theoremstyle{remark}

\newtheorem{remark}[theorem]{Remark}

\numberwithin{equation}{section}

\usepackage{color}

\usepackage[normalem]{ulem}

\def\Xint#1{\mathchoice
{\XXint\displaystyle\textstyle{#1}}%
{\XXint\textstyle\scriptstyle{#1}}%
{\XXint\scriptstyle
\scriptscriptstyle{#1}}%
{\XXint\scriptscriptstyle
\scriptscriptstyle{#1}}%
\!\int}
\def\XXint#1#2#3{{
\setbox0=\hbox{$#1{#2#3}{\int}$}
\vcenter{\hbox{$#2#3$}}\kern-.5\wd0}}
\def\dashint{\Xint-}

\begin{document}
\title{Convergence analysis of Hermite approximations for analytic functions}
\author{Haiyong Wang\footnotemark[1]~\footnotemark[2] ~ and ~ Lun Zhang\footnotemark[3]}
\date{}
\maketitle
\renewcommand{\thefootnote}{\fnsymbol{footnote}}

\footnotetext[1]{School of Mathematics and Statistics, Huazhong
University of Science and Technology, Wuhan 430074, P. R. China.
E-mail: \texttt{haiyongwang@hust.edu.cn}}

\footnotetext[2]{Hubei Key Laboratory of Engineering Modeling and
Scientific Computing, Huazhong University of Science and Technology,
Wuhan 430074, P. R. China}

\footnotetext[3]{School of Mathematical Sciences, Center for Applied Mathematics and Shanghai Key Laboratory for Contemporary Applied Mathematics, Fudan University, Shanghai 200433, P. R. China. E-mail:
\texttt{lunzhang@fudan.edu.cn}}

\begin{abstract}
In this paper, we present a rigorous analysis for root-exponential convergence of Hermite approximations, including projection and interpolation methods, for functions that are analytic in an infinite strip containing the real axis and satisfy certain restrictions on the asymptotic behavior at infinity within this strip. The key ingredients of our analysis are some new and remarkable contour integral representations for the Hermite coefficients and the remainder of Hermite spectral interpolations with which sharp error estimates for Hermite approximations in the weighted and maximum norms are established. Further extensions to Gauss-Hermite quadrature and the scaling factor are also discussed. Particularly, we prove the root-exponential convergence of Gauss--Hermite quadrature under explicit conditions on the integrands. Numerical experiments confirm our theoretical results.
\end{abstract}


\vspace{0.05in}


\section{Introduction}\label{sec:introduction}
Hermite polynomials/functions are orthogonal systems on the real line, which serve as a natural choice of basis functions when dealing with problems on unbounded domains. Hermite spectral methods, which are achieved by Hermite spectral approximations including expansions or interpolations via Hermite polynomials/functions, play an important role in a variety of physical applications such as quantum dynamics and plasma physics (cf. \cite{Boyd2000,Canuto2006,Shen2011}). Several attractive advantages of Hermite spectral methods can be summarized as follows:
\begin{itemize}
\item Hermite functions are the eigenfunctions or the limiting asymptotic eigenfunctions of many problems of physical interests \cite{Boyd1984,Lubich2008}.

\item Hermite differentiation matrix is skew-symmetric, tridiagonal and irreducible, which leads to conservative semidiscretisations for time-dependent PDEs \cite{Iserles2019}. Moreover, the spectral radii for the first and second Hermite differentiation matrices are $O(\sqrt{n})$ and $O(n)$, respectively, where $n$ is the number of terms of Hermite spectral approximation. This particularly implies rather weak stability restrictions on the time step when explicit time stepping methods are applied to the semidiscretisation system \cite{Weideman1992}.
\end{itemize}

A fundamental issue toward understanding the convergence behavior of Hermite spectral methods is to establish rigorous error analysis of Hermite spectral approximations. When the underlying function is analytic, Hille in \cite{Hille1939,Hille1940} and Szeg\H{o} in \cite{Szego1939} established the domain of convergence for the Hermite expansions in the complex plane. More precisely, Hille showed that the domain is an infinite strip when the underlying function is analytic and decays exponentially fast at infinity within this strip. Elliott and Tuan in \cite{Elliott1974} established a contour integral representation for the Hermite coefficients. Unfortunately, no further analysis on the convergence rate of Hermite approximations was carried out. Based on Hille's works, Boyd in \cite{Boyd1980,Boyd1984} refined the connection between the asymptotic estimate of Hermite coefficients and the behavior of the underlying functions at infinity. In particular, asymptotics of Hermite coefficients for some specific entire functions (namely, super- and sub-Gaussian types) and analytic functions with poles were obtained by using the steepest descent method. Although these studies provide important insights onto the convergence properties of Hermite approximations, a number of essential issues still remain open. For example, except for some special cases, no general convergence result was studied for Hermite projections and, as far as we know, no discussion was devoted to the convergence rate of Hermite spectral interpolations.

It is the aim of this work to present a rigorous analysis of Hermite spectral approximations, including Hermite projection and interpolation methods, for analytic functions. When the underlying function is analytic in an infinite strip containing the real axis and satisfies some restrictions on the asymptotic behavior at infinity, we establish contour integral representations for the Hermite coefficients and the remainder of Hermite spectral interpolations over the boundary of the strip. These, together with asymptotic behaviors of Hermite polynomials as well as their weighted Cauchy transforms, allow us to establish some sharp error estimates for Hermite projection and interpolation methods. Extensions of our analysis to Gauss-Hermite quadrature and the scaling factor of Hermite approximation are also discussed. Notably, we establish the root-exponential convergence of Gauss-Hermite quadrature under explicit conditions on the integrands, which has been considered as a folklore for a long time.

It is worthwhile to point out that rigorous analysis of spectral approximations have attracted renewed interests over the past few decades; cf. \cite{Reddy2005,Wang2012,Wang2021,Wang2023a,Wang2023c,Wang2014,Xiang2012a,Xiang2012b,Xie2013,Zhao2013}. Most of these studies, however, were devoted to Jacobi approximations as well as their special cases like Chebyshev, Legendre and Gegenbauer approximations. Compared with Jacobi case, the analysis of Laguerre and Hermite cases is much more difficult since their convergence rates depend not only on the locations of the singularities, but also on the asymptotic behavior of the underlying function at infinity. More recently, by exploiting contour integral techniques, the first-named author established some sharp error estimates for Laguerre projection and interpolation methods in \cite{Wang2023c}. Here we extend the idea therein to Hermite case and our analysis relies heavily on some remarkable contour integral representations for Hermite coefficients and the remainder of Hermite spectral interpolations. These integrals involve Hermite polynomials and their weighted Cauchy transforms, which might be of independent interests.

The rest of this paper is organized as follows. In the next section, we collect some properties of Hermite polynomials and present some important properties for their weighted Cauchy transform. In Sections \ref{sec:Projection} and \ref{sec:Interp}, we prove the root-exponential convergence of Hermite projections and interpolations for analytic functions, respectively. We extend our analysis to two topics of practical interests in Section \ref{sec:Extension} and present some concluding remarks in Section \ref{sec:Conclusion}.

\section{Hermite polynomials and their weighted Cauchy transforms}\label{sec:Cauchy}
In this section, we review some basic facts of Hermite
polynomials and establish some new properties of their weighted Cauchy transforms.
For more information about Hermite polynomials, we refer to
\cite{Olver2010,Szego1939} and the references therein.

Let $\omega(x)=\mathrm{e}^{-x^2}$ be the Hermite weight function and
$\mathbb{N}_0:=\{0,1,2,\ldots \}$. The Hermite polynomials are defined by
\begin{equation}\label{def:HermPoly}
H_n(x) = \frac{(-1)^n}{\omega(x)} \frac{\mathrm{d}^n\omega(x)}{\mathrm{d}x^n}, \quad n\in\mathbb{N}_0,
\end{equation}
and it is well-known that
\begin{equation}\label{eq:ortho}
\int_{-\infty}^{\infty} H_n(x) H_m(x) \omega(x) \mathrm{d}x = \gamma_{n} \delta_{n,m}, \quad n,m\in \mathbb{N}_0,
\end{equation}
where $\gamma_{n}= 2^n n!\sqrt{\pi} $ and $\delta_{n,m}$ stands for the Kronecker delta. Below we list some properties of Hermite polynomials that will be used later:
\begin{itemize}
  \item The three term recurrence relation of $H_n$ reads
  \begin{equation}\label{eq:recur}
  \begin{aligned}
  & H_{n+1}(x) = 2xH_n(x) - 2nH_{n-1}(x), \quad n\in\mathbb{N}, \\
  & H_0(x) = 1, \quad H_1(x) = 2x.
  \end{aligned}
  \end{equation}

  \item By \cite{Indritz1961}, we have for $x\in(-\infty,\infty)$
  \begin{equation}\label{eq:HermBound}
   \mathrm{e}^{-x^2/2} |H_n(x)| \leq \sqrt{2^nn!} = \pi^{-1/4} \sqrt{\gamma_{n}}.
  \end{equation}
  Moreover, by \cite[Theorem~8.22.7]{Szego1939}, we have for $z\in\mathbb{C}$ and $\Im(z)\neq0$
  \begin{equation}\label{eq:HermFunAsy}
  |\mathrm{e}^{-z^2/2}H_{n}(z)| = \frac{\Gamma(n+1)}{2\Gamma(n/2+1)}
  \mathrm{e}^{|\Im(z)|\sqrt{2n+1}} \left( 1 + O(n^{-1/2}) \right),
  \end{equation}
  where $\Gamma(z)$ is the gamma function, and it holds uniformly for bounded $z$.

  \item Since $\omega(x)=\omega(-x)$, it follows for $x\in\mathbb{R}$ that
  \begin{equation}\label{eq:HermSymm}
  H_n(-x)=(-1)^nH_n(x).
  \end{equation}
 Moreover, this property also holds for $z\in\mathbb{C}$, i.e., $H_n(-z)=(-1)^nH_n(z)$.
\end{itemize}

We then define the weighted Cauchy transform of Hermite polynomials by
\begin{equation}\label{def:Psi}
\Phi_n(z) = \frac{1}{2\pi \mathrm{i}}\int_{\mathbb{R}}
\frac{\omega(x) H_n(x)}{z-x} \mathrm{d}x, \quad
z\in\mathbb{C}\setminus\mathbb{R},
\end{equation}
where $\mathrm{i}$ is the imaginary unit. It is easily seen that $\Phi_n(z)$ is analytic in the whole complex plane with a cut along $\mathbb{R}$, and satisfies the same recurrence relation as $H_n(x)$ given in \eqref{eq:recur}. We establish some important properties of $\Phi_n(z)$ in what follows for later use.
\begin{lemma}\label{lem:Phin}
The function $\Phi_n$ defined in \eqref{def:Psi} satisfies the following properties.
\begin{itemize}\label{lem:Psi}
\item[\rm(i)] For any $x\in\mathbb{R}$ and each $n\in\mathbb{N}_0$, we have
\begin{equation}\label{eq:Plemelj}
\begin{cases}
{\displaystyle \lim_{\varepsilon\rightarrow0+} \left[
\Phi_n(x-\mathrm{i}\varepsilon) + \Phi_n(x+\mathrm{i}\varepsilon)
\right] = \frac{1}{\pi \mathrm{i}}\dashint_{\mathbb{R}}
\frac{\omega(y)H_n(y)}{x-y} \mathrm{d}y},  \\[3ex]
{\displaystyle \lim_{\varepsilon\rightarrow0+} \left[
\Phi_n(x-\mathrm{i}\varepsilon) - \Phi_n(x+\mathrm{i}\varepsilon)
\right] = \omega(x) H_n(x)},
\end{cases}
\end{equation}
where the bar indicates the Cauchy principal value.

\item[\rm(ii)] For each $n\in\mathbb{N}_0$, as $z\rightarrow\infty$ in $\epsilon < | \arg(z)| < \pi - \epsilon$ with $\epsilon \in (0, \pi)$, we have
\begin{equation}\label{eq:PhiAsyZ}
\Phi_n(z) = \frac{n!}{2\sqrt{\pi}\mathrm{i}}z^{-n-1} + O(z^{-n-3}).
\end{equation}

\item[\rm(iii)] Let $U(a,z)$ be the Weber parabolic cylinder function \cite[Chapter~12]{Olver2010}. For each $n\in\mathbb{N}_0$, we have
\begin{align}\label{eq:PhiFormula}
\Phi_n(z) &= \frac{2^{n/2} n! \mathrm{e}^{-z^2/2}}{\sqrt{2\pi}} \left\{
\begin{array}{ll}
{\displaystyle U\left(n+\frac{1}{2},-\mathrm{i}\sqrt{2}z\right) \mathrm{e}^{-\mathrm{i}(n+2)\pi/2} },    & \Im(z)>0,  \\[3ex]
{\displaystyle U\left(n+\frac{1}{2},\mathrm{i}\sqrt{2}z\right) \mathrm{e}^{\mathrm{i}n\pi/2} },         & \Im(z)<0.
\end{array}
\right.
\end{align}
Moreover, let $U(a,b,z)$ be the confluent hypergeometric function of the second kind (also known as Kummer's function of the second kind, cf. \cite[Chapter~13]{Olver2010}), we also have
\begin{equation}\label{eq:PhiFormula1}
\Phi_n(z) = \frac{n!}{2\sqrt{\pi}}
U\left(\frac{n+1}{2},\frac{1}{2}, -z^2\right) \left\{
\begin{array}{ll}
{\displaystyle \mathrm{e}^{-\mathrm{i}(n+2)\pi/2}},    & \Im(z)>0,  \\[2ex]
{\displaystyle \mathrm{e}^{\mathrm{i}n\pi/2}},         & \Im(z)<0.
\end{array}
\right.
\end{equation}

\item[\rm(iv)] As $n\rightarrow\infty$, we have
\begin{align}\label{eq:PhiAsy}
\Phi_n(z) &\sim \frac{\Gamma(2\kappa)}{2\Gamma(\kappa)\sqrt{\kappa}} \exp\left(-\frac{z^2}{2} \pm 2\mathrm{i}\sqrt{\kappa}z \right) \left(\sum_{k=0}^{\infty} \frac{\Pi_k(z)}{\kappa^{k/2}} \right) \left\{
\begin{array}{ll}
{\displaystyle \mathrm{e}^{-\mathrm{i}(n+2)\pi/2}},    & \Im(z)>0,  \\[2ex]
{\displaystyle \mathrm{e}^{\mathrm{i}n\pi/2}},         & \Im(z)<0,
\end{array}
\right.
\end{align}
uniformly for bounded $z$, where $\kappa=(n+1)/2$, the $+$ and $-$ signs are taken for $\Im(z)>0$ and $\Im(z)<0$, respectively. The functions $\Pi_k(z)$, $k\in\mathbb{N}_0$, are polynomials in $z$ with
\[
\Pi_0(z)=1, \qquad \Pi_1(z) = \mathrm{i}\frac{z(z^2+3)}{12} \left\{
\begin{array}{ll}
{\displaystyle -1},    & \Im(z)>0,  \\[2ex]
{\displaystyle 1},     & \Im(z)<0.
\end{array}
\right.
\]
\end{itemize}
\end{lemma}
\begin{proof}
The first statement is a direct consequence of the Plemelj formula \cite{Muskhl1972}. For the second statement, by expanding $1/(z-x)$ in \eqref{def:Psi} as
\[
\frac{1}{z-x} = \sum_{k=0}^{n+1} \frac{x^k}{z^{k+1}} + \frac{x^{n+2}}{(z-x)z^{n+2}},
\]
and using the orthogonality condition \eqref{eq:ortho} and the symmetry \eqref{eq:HermSymm}, it is readily seen that
\begin{align}
\Phi_n(z) &= \frac{z^{-n-1}}{2\pi\mathrm{i}}\int_{\mathbb{R}} x^n H_n(x) \omega(x) \mathrm{d}x + \frac{z^{-n-2}}{2\pi\mathrm{i}} \int_{\mathbb{R}} \frac{x^{n+2} H_n(x) \omega(x)}{z-x} \mathrm{d}x \nonumber \\
&= \frac{n!}{2\sqrt{\pi}\mathrm{i}} z^{-n-1} + O(z^{-n-3}), \nonumber
\end{align}
where we have used the fact that the leading coefficient of $H_n(x)$ is $2^n$ in the last step. This proves \eqref{eq:PhiAsyZ}. Alternatively, assuming \eqref{eq:PhiFormula1}, it also follows from the fact (cf. \cite[Equation (13.7.3)]{Olver2010}) that as $z\to\infty$,
$$
U(a,b,z)=z^{-a}\left(1+O(z^{-1}) \right), \qquad |\arg(z)| < \frac{3\pi}{2}.
$$
For the third statement, by the integral representation of Hermite polynomials given in \cite[Equation~(18.10.10)]{Olver2010} and exchanging the order of integration, we have
\begin{align}\label{eq:PhiPara}
\Phi_n(z) &= \frac{2^{n+1}}{2\pi\mathrm{i}\sqrt{\pi}} \int_{0}^{\infty} \mathrm{e}^{-t^2}t^n \left( \int_{\mathbb{R}} \frac{\cos(2xt - n\pi/2)}{z-x}\mathrm{d}x \right) \mathrm{d}t  \nonumber \\[1ex]
&= \frac{2^{n}}{\sqrt{\pi}} \left\{
\begin{array}{ll}
{\displaystyle \mathrm{e}^{-(n+2)\pi\mathrm{i}/2} \int_{0}^{\infty} \mathrm{e}^{-t^2 + 2\mathrm{i}zt} t^n \mathrm{d}t },    & \Im(z)>0,  \\[2ex]
{\displaystyle \mathrm{e}^{n\pi\mathrm{i}/2} \int_{0}^{\infty} \mathrm{e}^{-t^2-2\mathrm{i}zt} t^n \mathrm{d}t },           & \Im(z)<0,
\end{array}
\right. \nonumber
\\
&= \frac{2^{n/2}}{\sqrt{2\pi} } \left\{
\begin{array}{ll}
{\displaystyle \mathrm{e}^{-(n+2)\pi\mathrm{i}/2} \int_{0}^{\infty} \mathrm{e}^{-t^2/2 + \mathrm{i}\sqrt{2}zt} t^n \mathrm{d}t },    & \Im(z)>0,  \\[2ex]
{\displaystyle \mathrm{e}^{n\pi\mathrm{i}/2} \int_{0}^{\infty} \mathrm{e}^{-t^2/2-\mathrm{i}\sqrt{2}zt} t^n \mathrm{d}t },           & \Im(z)<0,
\end{array}
\right.
\end{align}
where the second equality follows from the facts that
\[
\int_{\mathbb{R}} \frac{\mathrm{e}^{\mathrm{i}\lambda x}}{z-x} \mathrm{d}x = \left\{
\begin{array}{ll}
{\displaystyle -2\pi \mathrm{i} \mathrm{e}^{\mathrm{i}\lambda z}},    & \Im(z)>0,  \\[2ex]
{\displaystyle 0 },         & \Im(z)<0,
\end{array}
\right. \quad \int_{\mathbb{R}} \frac{\mathrm{e}^{-\mathrm{i}\lambda x}}{z-x} \mathrm{d}x = \left\{
\begin{array}{ll}
{\displaystyle 0},    & \Im(z)>0,  \\[2ex]
{\displaystyle 2\pi \mathrm{i} \mathrm{e}^{-\mathrm{i}\lambda z}},         & \Im(z)<0,
\end{array}
\right.
\]
for $\lambda>0$, which can be obtained by the residue theorem and Jordan's lemma, and the last equality follows from the second one by a change of variable $t\mapsto t/\sqrt{2}$. The equality \eqref{eq:PhiFormula} then follows by combining the last equality  of \eqref{eq:PhiPara} with the integral representation of $U(a,z)$ \cite[Equation~(12.5.1)]{Olver2010}. A combination of \eqref{eq:PhiFormula} and \cite[Equation~(12.7.14)]{Olver2010} gives us \eqref{eq:PhiFormula1}. 

Finally, to show the asymptotic expansion \eqref{eq:PhiAsy}, by \cite[Equation~(10.3.37)]{TemmeBook2015} we know that for $z\in \mathbb{C}\setminus (-\infty,0]$,
\begin{equation}
U(a,b,z) \sim 2\left(\frac{z}{a}\right)^{(1-b)/2} \frac{\mathrm{e}^{z/2}}{\Gamma(a)} \sum_{k=0}^{\infty} c_k(z) \left(\frac{z}{a}\right)^{k/2} K_{k+1-b}(2\sqrt{az}), \quad a\to +\infty, \nonumber
\end{equation}
uniformly for bounded $z$ and $b$, where $K_{\nu}(z)$ is the modified Bessel function of the second kind and $c_k(z)$, $k\in\mathbb{N}_0$, are polynomials in $z$ of degree $k$ whose first few terms can be found in \cite[Equation~(10.3.31)]{TemmeBook2015}. 
On account of the fact (cf. \cite[Equation (10.40.2)]{Olver2010}) that as $z\rightarrow\infty$
$$
K_{\nu}(z) \sim \mathrm{e}^{-z} \sqrt{\frac{\pi}{2z}} \sum_{k=0}^{\infty}\frac{\phi_k(\nu)}{z^k}, \qquad |\arg(z)|<\frac{3\pi}{2},
$$
where $\phi_0(\nu)=1$ and $\phi_k(\nu)=\prod_{j=1}^{k}(4\nu^2-(2j-1)^2)/(8^kk!)$ for $k\geq1$, it is readily verified that as $a\rightarrow+\infty$,
\begin{align}\label{eq:AsyUFull}
U(a,b,z) &\sim \frac{\sqrt{\pi}}{\Gamma(a)}\left(\frac{z}{a}\right)^{(1-b)/2}\frac{\mathrm{e}^{z/2-2\sqrt{az}}}{(az)^{1/4}}
\sum_{k=0}^{\infty} c_k(z) \left(\frac{z}{a}\right)^{k/2} \sum_{j=0}^{\infty} \frac{\phi_{j}(k+1-b)}{(2\sqrt{az})^{j}}.
\end{align}
Combining this expansion with \eqref{eq:PhiFormula1} and noting $\phi_j(k+1/2)=0$ for $j\geq k+1$ gives us \eqref{eq:PhiAsy}. This completes the proof of Lemma \ref{lem:Phin}.
\end{proof}

\begin{remark}
By keeping the leading term in \eqref{eq:AsyUFull}, it is readily seen that
\begin{equation}\label{eq:AsyU}
U(a,b,z) = \frac{\sqrt{\pi}}{\Gamma(a)}\left(\frac{z}{a}\right)^{(1-b)/2}\frac{\mathrm{e}^{z/2-2\sqrt{az}}}{(az)^{1/4}}
\left(1 + O(a^{-1/2})  \right), \quad a\to +\infty,
\end{equation}
uniformly for bounded $z$ and $b$. As this approximation will be used frequently in subsequent sections, we record it here.
\end{remark}

\begin{remark}\label{rk:boundary}
For any $x\in\mathbb{R}$, from \eqref{eq:PhiFormula} we obtain
\begin{align}
\lim_{\varepsilon\rightarrow0+} \Phi_n(x\pm\mathrm{i}\varepsilon) = \frac{2^{n/2} n! \mathrm{e}^{-x^2/2}}{\sqrt{2\pi}} U\left(n+\frac{1}{2},\mp\mathrm{i}\sqrt{2}x\right) \left(\mp\mathrm{e}^{\mp\mathrm{i}n\pi/2}\right).  \nonumber
\end{align}
As $x\rightarrow\pm\infty$, by the asymptotic of $U(a,z)$ as $z\rightarrow\infty$ in \cite[Equation~(12.9.1)]{Olver2010}, it follows that
\begin{align}
\lim_{\varepsilon\rightarrow0+} \Phi_n(x\pm\mathrm{i}\varepsilon) &= \left(\mp\mathrm{e}^{\mp\mathrm{i}n\pi/2}\right) \frac{2^{n/2} n!}{\sqrt{2\pi}} \left(\mp\mathrm{i}\sqrt{2}x\right)^{-n-1} \left(1 + O(x^{-2})\right)  \nonumber \\
&= \frac{n!}{2\sqrt{\pi}\mathrm{i}} x^{-n-1} \left(1 + O(x^{-2})\right).  \nonumber
\end{align}
\end{remark}

\begin{remark}
Using a different approach, Rusev also studied the large $n$ behavior of $\Phi_n(z)$ in \cite{Rusev}. The difference between the result therein and \eqref{eq:PhiFormula1} lies in that Rusev's asymptotic results involve the term $1+k_n(z)$ with $k_n(z)$ being complex functions analytic in the upper/lower half plane and $\lim_{n\rightarrow\infty}k_n(z)=0$ uniformly on every compact set of the upper/lower half plane, respectively. Our result gives a more precise expansion for this term, i.e., $\sum_{k=0}^{\infty} \Pi_k(z)/\kappa^{k/2}$.
\end{remark}

\section{Convergence analysis of Hermite projections}\label{sec:Projection}
In this section, we carry out a thorough convergence analysis of Hermite projections. To proceed, we introduce an infinite strip in the complex plane
\begin{equation}\label{def:Strip}
\mathcal{S}_{\rho} := \big\{z\in\mathbb{C}: ~\Im(z)\in[-\rho,\rho]
\big\},
\end{equation}
where $\rho>0$, and denote by $\partial{\mathcal{S}}_{\rho}$ the boundary of $\mathcal{S}_{\rho}$, i.e., $\partial{\mathcal{S}}_{\rho}=\{z\in\mathbb{C}:~z=x\pm\mathrm{i}\rho, ~x\in(-\infty,\infty)\}$; see Figure \ref{fig:strip} for an illustration. A key motivation for introducing $\mathcal{S}_{\rho}$ is that the domain of convergence of Hermite approximations is actually characterized by $\mathcal{S}_{\rho}$; see \cite{Hille1940} and \cite[Chapter~9]{Szego1939}.

\begin{figure}[t]
\centering

\tikzset{every picture/.style={line width=0.75pt}} 

\begin{tikzpicture}[x=0.75pt,y=0.75pt,yscale=-1,xscale=1]

\draw  [fill={rgb, 255:red, 155; green, 155; blue, 155 }  ,fill opacity=1 ][dash pattern={on 0.84pt off 2.51pt}] (114.58,66.88) -- (294.55,66.88) -- (294.55,202.09) -- (114.58,202.09) -- cycle ;
\draw    (81.75,134.87) -- (325.38,134.12) ;
\draw [shift={(327.38,134.11)}, rotate = 179.82] [fill={rgb, 255:red, 0; green, 0; blue, 0 }  ][line width=0.08]  [draw opacity=0] (12,-3) -- (0,0) -- (12,3) -- cycle    ;
\draw  [dash pattern={on 4.5pt off 4.5pt}]  (115.02,202.09) -- (294.55,202.09) ;
\draw [shift={(211.79,202.09)}, rotate = 180] [fill={rgb, 255:red, 0; green, 0; blue, 0 }  ][line width=0.08]  [draw opacity=0] (12,-3) -- (0,0) -- (12,3) -- cycle    ;
\draw  [dash pattern={on 4.5pt off 4.5pt}]  (114.58,66.88) -- (294.11,66.88) ;
\draw [shift={(195.35,66.88)}, rotate = 0] [fill={rgb, 255:red, 0; green, 0; blue, 0 }  ][line width=0.08]  [draw opacity=0] (12,-3) -- (0,0) -- (12,3) -- cycle    ;
\draw    (203.97,24.21) -- (203.97,247.33) ;
\draw [shift={(203.97,22.21)}, rotate = 90] [fill={rgb, 255:red, 0; green, 0; blue, 0 }  ][line width=0.08]  [draw opacity=0] (12,-3) -- (0,0) -- (12,3) -- cycle    ;

\draw (171,24.17) node [anchor=north west][inner sep=0.75pt]   [align=left] {$\Im (z)$};
\draw (306.53,137.55) node [anchor=north west][inner sep=0.75pt]   [align=left] {$\Re (z)$};
\draw (234.53,87.55) node [anchor=north west][inner sep=0.75pt]   [align=left] {$\mathcal{S}_{\rho}$};

\draw (234.53,203.55) node [anchor=north west][inner sep=0.75pt]   [align=left] {$\partial \mathcal{S}_{\rho}$};

\end{tikzpicture}

\caption{The infinite strip $\mathcal{S}_{\rho}$ and its boundary $\partial \mathcal{S}_{\rho}$.}
\label{fig:strip}
\end{figure}

\subsection{Projections using Hermite polynomials}
Let $\mathbb{P}_n$ be the space of polynomials of degree up to $n$ and denote by $L_{\omega}^2(\mathbb{R})$ the Hilbert space consisting of all squared integrable functions with respect to the Hermite weight function $\omega(x)=\mathrm{e}^{-x^2}$ over $\mathbb{R}$, equipped with the inner product and the norm
\begin{equation}\label{def:Lomega2}
\langle f,g \rangle_{\omega} := \int_{-\infty}^{\infty} f(x) g(x) \omega(x) \mathrm{d}x, \quad ||f||_{L_{\omega}^2(\mathbb{R})}=\sqrt{\langle f,f \rangle_{\omega}}.
\end{equation}
We set $\Pi_n^{\mathrm{P}}:L_{\omega}^2(\mathbb{R})\rightarrow\mathbb{P}_n$ to be the orthogonal projection operator, i.e., for any $f\in{L}_{\omega}^2(\mathbb{R})$,
\begin{equation}\label{eq:HermiteExp}
(\Pi_n^\mathrm{P}f)(x) = \sum_{k=0}^{n} a_k H_k(x), \qquad  a_k =
\frac{\langle f, H_k \rangle_{\omega}}{\gamma_{k}},
\end{equation}
where $H_k(x)$ is the Hermite polynomial of degree $k$. The following theorem shows that the coefficient $a_n$ admits a contour integral representation involving the weighted Cauchy transform of the Hermite polynomial $H_n(x)$, which plays an important role in our further analysis. In what follows, the orientation of a closed curve in the complex plane is always taken in a counterclockwise manner.

\begin{theorem}\label{thm:HermCoeff}
If $f$ is analytic in the strip $\mathcal{S}_{\rho}$\footnote{It means that $f$ is analytic in an open strip that contains $\mathcal{S}_{\rho}$, thus including the boundary $\Im(z)=\pm\rho$. This assumption might be relaxed to the one that $f$ is analytic in the strip $\mathcal{S}_{\rho-\epsilon}$ for arbitrarily small $\epsilon>0$ and continuous up to the boundary $\Im(z)=\pm\rho$, at the cost of more cumbersome proofs.} for
some $\rho>0$ and $|f(z)|\leq\mathcal{K}|z|^{\sigma}$ for some
$\sigma\in\mathbb{R}$ as $|z|\rightarrow\infty$ within the strip, it then holds that, for each
$n\geq\max\{\lfloor\sigma\rfloor+1,0\}$,
\begin{equation}\label{eq:HermCoeff}
a_n = \frac{1}{\gamma_n} \int_{\partial{\mathcal{S}}_{\rho}}
\Phi_n(z) f(z) \mathrm{d}z,
\end{equation}
where $\Phi_n$ is defined in \eqref{def:Psi}.
\end{theorem}
\begin{proof}
Let $\varepsilon$ and $\eta$ be two arbitrary positive constants, we set the following
auxiliary contours:
\begin{equation}\label{def:contours}
\begin{aligned}
\Gamma_{\pm}^{\mathrm{D}}:~~ z &= x \pm \mathrm{i}\varepsilon, \quad \Gamma_{\pm}^{\mathrm{U}}:~~ z = x \pm \mathrm{i}\rho, \qquad
x\in[-\eta,\eta],
 \\[2ex]
\Gamma_{\pm}^{\mathrm{R}}:~~ z & = \eta \pm
\mathrm{i}y,
 \quad \Gamma_{\pm}^{\mathrm{L}}:~~ z = -\eta \pm
\mathrm{i}y, \qquad y\in[\varepsilon,\rho];
\end{aligned}
\end{equation}
see Figure \ref{figure-Gamma} for an illustration and the orientations.

\begin{figure}[t]
\centering
\tikzset{every picture/.style={line width=0.75pt}} 
\begin{tikzpicture}[x=0.75pt,y=0.75pt,yscale=-1.2,xscale=1.2]
\draw    (167,120.33) -- (369.89,120.77) ;
\draw [shift={(371.89,120.78)}, rotate = 180.12] [fill={rgb, 255:red, 0; green, 0; blue, 0 }  ][line width=0.08]  [draw opacity=0] (12,-3) -- (0,0) -- (12,3) -- cycle    ;
\draw    (190.78,200.46) -- (190.23,32.35) ;
\draw [shift={(190.22,30.35)}, rotate = 89.81] [fill={rgb, 255:red, 0; green, 0; blue, 0 }  ][line width=0.08]  [draw opacity=0] (12,-3) -- (0,0) -- (12,3) -- cycle    ;
\draw    (229.44,70.76) -- (330.11,70.78) ;
\draw [shift={(270.78,70.77)}, rotate = 0.01] [fill={rgb, 255:red, 0; green, 0; blue, 0 }  ][line width=0.08]  [draw opacity=0] (12,-3) -- (0,0) -- (12,3) -- cycle    ;
\draw    (229.44,70.76) -- (229.44,100.11) ;
\draw [shift={(229.44,92.44)}, rotate = 270] [fill={rgb, 255:red, 0; green, 0; blue, 0 }  ][line width=0.08]  [draw opacity=0] (12,-3) -- (0,0) -- (12,3) -- cycle    ;
\draw    (229.44,100.11) -- (330.11,100.12) ;
\draw [shift={(286.78,100.12)}, rotate = 180.01] [fill={rgb, 255:red, 0; green, 0; blue, 0 }  ][line width=0.08]  [draw opacity=0] (12,-3) -- (0,0) -- (12,3) -- cycle    ;
\draw    (330.11,70.78) -- (330.11,100.12) ;
\draw [shift={(330.11,76.45)}, rotate = 90] [fill={rgb, 255:red, 0; green, 0; blue, 0 }  ][line width=0.08]  [draw opacity=0] (12,-3) -- (0,0) -- (12,3) -- cycle    ;
\draw   (189.78,69.78) .. controls (189.78,70.15) and (190.08,70.44) .. (190.44,70.44) .. controls (190.81,70.44) and (191.11,70.15) .. (191.11,69.78) .. controls (191.11,69.41) and (190.81,69.11) .. (190.44,69.11) .. controls (190.08,69.11) and (189.78,69.41) .. (189.78,69.78) -- cycle ;
\draw   (189.78,100.11) .. controls (189.78,100.48) and (190.08,100.78) .. (190.44,100.78) .. controls (190.81,100.78) and (191.11,100.48) .. (191.11,100.11) .. controls (191.11,99.74) and (190.81,99.44) .. (190.44,99.44) .. controls (190.08,99.44) and (189.78,99.74) .. (189.78,100.11) -- cycle ;
\draw   (229.78,120.44) .. controls (229.78,120.81) and (230.08,121.11) .. (230.44,121.11) .. controls (230.81,121.11) and (231.11,120.81) .. (231.11,120.44) .. controls (231.11,120.08) and (230.81,119.78) .. (230.44,119.78) .. controls (230.08,119.78) and (229.78,120.08) .. (229.78,120.44) -- cycle ;
\draw   (329.78,120.94) .. controls (329.78,121.31) and (330.08,121.61) .. (330.44,121.61) .. controls (330.81,121.61) and (331.11,121.31) .. (331.11,120.94) .. controls (331.11,120.58) and (330.81,120.28) .. (330.44,120.28) .. controls (330.08,120.28) and (329.78,120.58) .. (329.78,120.94) -- cycle ;
\draw    (230.24,140.64) -- (330.91,140.65) ;
\draw [shift={(271.58,140.65)}, rotate = 0.01] [fill={rgb, 255:red, 0; green, 0; blue, 0 }  ][line width=0.08]  [draw opacity=0] (12,-3) -- (0,0) -- (12,3) -- cycle    ;
\draw    (230.24,140.64) -- (230.24,169.99) ;
\draw [shift={(230.24,162.31)}, rotate = 270] [fill={rgb, 255:red, 0; green, 0; blue, 0 }  ][line width=0.08]  [draw opacity=0] (12,-3) -- (0,0) -- (12,3) -- cycle    ;
\draw    (230.24,169.99) -- (330.91,170) ;
\draw [shift={(287.58,170)}, rotate = 180.01] [fill={rgb, 255:red, 0; green, 0; blue, 0 }  ][line width=0.08]  [draw opacity=0] (12,-3) -- (0,0) -- (12,3) -- cycle    ;
\draw    (330.91,140.65) -- (330.91,170) ;
\draw [shift={(330.91,146.33)}, rotate = 90] [fill={rgb, 255:red, 0; green, 0; blue, 0 }  ][line width=0.08]  [draw opacity=0] (12,-3) -- (0,0) -- (12,3) -- cycle    ;
\draw   (190.04,139.93) .. controls (190.04,140.21) and (190.26,140.43) .. (190.54,140.43) .. controls (190.81,140.43) and (191.03,140.21) .. (191.03,139.93) .. controls (191.03,139.66) and (190.81,139.44) .. (190.54,139.44) .. controls (190.26,139.44) and (190.04,139.66) .. (190.04,139.93) -- cycle ;
\draw   (189.76,170.16) .. controls (189.76,170.53) and (190.06,170.83) .. (190.43,170.83) .. controls (190.79,170.83) and (191.09,170.53) .. (191.09,170.16) .. controls (191.09,169.79) and (190.79,169.49) .. (190.43,169.49) .. controls (190.06,169.49) and (189.76,169.79) .. (189.76,170.16) -- cycle ;
\draw (178,64) node [anchor=north west][inner sep=0.75pt]   [align=left] {$\rho$};
\draw (179,95) node [anchor=north west][inner sep=0.75pt]   [align=left] {$\epsilon$};
\draw (267,50) node [anchor=north west][inner sep=0.75pt]   [align=left] {$\Gamma_+^{\mathrm{U}}$};
\draw (210,76) node [anchor=north west][inner sep=0.75pt]   [align=left] {$\Gamma_+^{\mathrm{L}}$};
\draw (333,76) node [anchor=north west][inner sep=0.75pt]   [align=left] {$\Gamma_+^{\mathrm{R}}$};
\draw (267,102) node [anchor=north west][inner sep=0.75pt]   [align=left] {$\Gamma_+^{\mathrm{D}}$};
\draw (267,123) node [anchor=north west][inner sep=0.75pt]   [align=left] {$\Gamma_-^{\mathrm{D}}$};
\draw (210,146.05) node [anchor=north west][inner sep=0.75pt]   [align=left] {$\Gamma_-^{\mathrm{L}}$};
\draw (333,146.05) node [anchor=north west][inner sep=0.75pt]   [align=left] {$\Gamma_-^{\mathrm{R}}$};
\draw (267,172) node [anchor=north west][inner sep=0.75pt]   [align=left] {$\Gamma_-^{\mathrm{U}}$};
\draw (219.34,120.05) node [anchor=north west][inner sep=0.75pt]   [align=left] {$-\eta$};
\draw (324,122.8) node [anchor=north west][inner sep=0.75pt]   [align=left] {$\eta$};
\draw (169,133) node [anchor=north west][inner sep=0.75pt]   [align=left] {$-\epsilon$};
\draw (167,163) node [anchor=north west][inner sep=0.75pt]   [align=left] {$-\rho$};
\end{tikzpicture}
\caption{The contours $\Gamma_{\pm}^{\mathrm{D}},\Gamma_{\pm}^{\mathrm{U}},\Gamma_{\pm}^{\mathrm{R}}$ and $\Gamma_{\pm}^{\mathrm{L}}$ defined in \eqref{def:contours}.}
\label{figure-Gamma}
\end{figure}

Since $\Phi_n(z)$ is analytic in the whole complex plane except the real line, it follows from Cauchy's theorem that
\begin{align}\label{eq:Cauchy}
\int_{\Gamma_{+}^{\mathrm{D}}\cup\Gamma_{-}^{\mathrm{D}}\cup\Gamma_{-}^{\mathrm{R}}\cup\Gamma_{+}^{\mathrm{R}}
\cup\Gamma_{+}^{\mathrm{U}}\cup\Gamma_{-}^{\mathrm{U}} \cup \Gamma_{+}^{\mathrm{L}} \cup \Gamma_{-}^{\mathrm{L}}}
\Phi_n(z) f(z) \mathrm{d}z &= 0.
\end{align}
We next evaluate the integrals over different contours as  $\varepsilon\rightarrow0^{+}$ and $\eta \to \infty$.
For the contour integrals over $\Gamma_{-}^{\mathrm{D}}\cup\Gamma_{+}^{\mathrm{D}}$, by Lemma
\ref{lem:Phin}, Remark \ref{rk:boundary} and the condition $|f(z)|\leq \mathcal{K}|z|^{\sigma}$ with $\sigma<\lfloor\sigma\rfloor+1\leq n$ for large $z$, we see for $x\in\mathbb{R}$ that
\[
|\Phi_n(x\pm\mathrm{i}\varepsilon)f(x\pm\mathrm{i}\varepsilon)| \leq \mathcal{K} (1 + |x\pm\mathrm{i}\varepsilon|)^{\sigma-n-1} \leq \mathcal{K} (1 + |x|)^{\sigma-n-1},
\]
where $\mathcal{K}$ is a generic positive constant and the last bound holds uniformly for small $\varepsilon$. Moreover, note that
\[
\int_{\mathbb{R}} (1 + |x|)^{\sigma-n-1} \mathrm{d}x = \frac{2}{n-\sigma} < \infty,
\]
we have, by the dominated convergence theorem, that
\begin{align}
\lim_{\substack{\varepsilon\rightarrow0^{+} \\
\eta\rightarrow\infty}} \int_{\Gamma_{-}^{\mathrm{D}}\cup\Gamma_{+}^{\mathrm{D}}}
\Phi_n(z) f(z) \mathrm{d}z &= \lim_{\substack{\varepsilon\rightarrow0^{+} \\
\eta\rightarrow\infty}} \int_{-\eta}^{\eta} \big( \Phi_n(x+\mathrm{i}\varepsilon)f(x+\mathrm{i}\varepsilon) - \Phi_n(x-\mathrm{i}\varepsilon)f(x-\mathrm{i}\varepsilon) \big) \mathrm{d}x \nonumber \\
&= \lim_{\eta\rightarrow\infty} \int_{-\eta}^{\eta} \lim_{\varepsilon\rightarrow0^{+}} \big( \Phi_n(x+\mathrm{i}\varepsilon)f(x+\mathrm{i}\varepsilon) - \Phi_n(x-\mathrm{i}\varepsilon)f(x-\mathrm{i}\varepsilon) \big) \mathrm{d}x \nonumber \\
&= \lim_{\eta\rightarrow\infty} \int_{-\eta}^{\eta} \lim_{\varepsilon\rightarrow0^{+}} \big( \Phi_n(x+\mathrm{i}\varepsilon) - \Phi_n(x-\mathrm{i}\varepsilon) \big) f(x) \mathrm{d}x \nonumber \\
&= - \lim_{\eta\rightarrow\infty} \int_{-\eta}^{\eta} \omega(x) H_n(x) f(x) \mathrm{d}x \nonumber \\
&= -\gamma_n a_n, \nonumber
\end{align}
where we have used the second equation of \eqref{eq:Plemelj} in the fourth equality and the definition of $a_n$ in \eqref{eq:HermiteExp} in the last equality. For the contour integral over the vertical line $\Gamma_{+}^{\mathrm{R}}$, by Lemma
\ref{lem:Phin}, Remark \ref{rk:boundary} and the condition $|f(z)|\leq \mathcal{K}|z|^{\sigma}$ with $\sigma<\lfloor\sigma\rfloor+1\leq n$ for large $z$ again, we see that for large $\eta$
\[
|\Phi_n(\eta+\mathrm{i}y) f(\eta+\mathrm{i}y)| \leq \mathcal{K} |\eta+\mathrm{i}y|^{\sigma-n-1} = \mathcal{K} (\eta^2 + y^2)^{(\sigma-n-1)/2},
\]
and by the dominated convergence theorem again, we have
\begin{align*}
\lim_{\substack{\varepsilon\rightarrow0^{+} \\
\eta\rightarrow\infty}}
\int_{\Gamma_{+}^{\mathrm{R}}} \Phi_n(z)
f(z) \mathrm{d}z &= \mathrm{i}
\lim_{\eta\rightarrow\infty}\int_{0^+}^{\rho}
\Phi_n(\eta+\mathrm{i}y) f(\eta+\mathrm{i}y) \mathrm{d}y = 0.
\end{align*}
Similarly, one has
\begin{align*}
\lim_{\substack{\varepsilon\rightarrow0^{+} \\
\eta\rightarrow\infty}}
\int_{\Gamma_{-}^{\mathrm{R}}} \Phi_n(z)
f(z) \mathrm{d}z =\lim_{\substack{\varepsilon\rightarrow0^{+} \\
\eta\rightarrow\infty}}
\int_{\Gamma_{+}^{\mathrm{L}}} \Phi_n(z)
f(z) \mathrm{d}z = \lim_{\substack{\varepsilon\rightarrow0^{+} \\
\eta\rightarrow\infty}}
\int_{\Gamma_{-}^{\mathrm{L}}} \Phi_n(z)
f(z) \mathrm{d}z=0.
\end{align*}
Finally, it is easily seen that
\begin{align}
\lim_{\eta\rightarrow\infty}
\int_{\Gamma_{+}^{\mathrm{U}}\cup \Gamma_{-}^{\mathrm{U}}} \Phi_n(z)
f(z) \mathrm{d}z &= \int_{\partial\mathcal{S}_{\rho}} \Phi_n(z) f(z)
\mathrm{d}z. \nonumber
\end{align}
Again, by item (ii) of Lemma \ref{lem:Phin} and the fact that $|f(z)|\leq \mathcal{K}|z|^{\sigma}$ with $\sigma<n$ as $|z|\rightarrow\infty$, the integral on the right-hand side above converges. By taking $\varepsilon\rightarrow0^{+}$ and $\eta \to \infty$ on both sides of \eqref{eq:Cauchy}, we arrive at \eqref{eq:HermCoeff}. This ends the proof of Theorem \ref{thm:HermCoeff}.
\end{proof}

\begin{remark}
It is easily seen that \eqref{eq:HermCoeff} holds for all $n\geq0$
when $|f(z)|\leq \mathcal{K}|z|^{\sigma}$ for some $\sigma<0$.
\end{remark}

As a consequence of the above theorem, we are able to establish an upper bound for the Hermite coefficient $a_n$ and an error bound for $\Pi_n^\mathrm{P}$ in the weighted $L^2$-norm. We start with an useful property of $\Phi_n(z)$.
\begin{lemma}\label{lem:MaxPhi}
For any $\rho>0$, we have
\begin{equation}\label{eq:PhiMax}
\max_{z\in\partial\mathcal{S}_{\rho}} |\Phi_n(z)| = |\Phi_n(\pm\mathrm{i}\rho)|.
\end{equation}
\end{lemma}
\begin{proof}
By \eqref{eq:HermSymm}, it is easily checked that $\Phi_n(z)=(-1)^{n+1}\Phi_n(-z)$, $z\in \mathbb{C} \setminus \mathbb{R}$. Thus, it suffices to show $\max_{\Im(z)=\rho} |\Phi_n(z)| = |\Phi_n(\mathrm{i}\rho)|$. Let $z=x+\mathrm{i}\rho$ with $x\in(-\infty,\infty)$, we see from \eqref{eq:PhiFormula} that
\begin{align}
|\Phi_n(z)| &= \frac{2^{n/2} n! \mathrm{e}^{-(x^2-\rho^2)/2}}{\sqrt{2\pi}} \left| U\left(n+\frac{1}{2},\sqrt{2}(\rho-\mathrm{i}x)\right) \right|. \nonumber
\end{align}
On the other hand, by the integral representation of $U(a,z)$ in \cite[Equation~(12.5.1)]{Olver2010}, we have that for $a>-1/2$ and $z\in\mathbb{C}$,
\begin{align}
|U(a,z)| &= \frac{|\mathrm{e}^{-z^2/4}|}{\Gamma(a+1/2)} \left| \int_{0}^{\infty} t^{a-1/2} \mathrm{e}^{-t^2/2-zt} \mathrm{d}t \right| \nonumber \\
&\leq \frac{\mathrm{e}^{-(\Re(z)^2-\Im(z)^2)/4}}{\Gamma(a+1/2)} \int_{0}^{\infty} t^{a-1/2} \mathrm{e}^{-t^2/2-\Re(z)t} \mathrm{d}t \nonumber \\
&= \mathrm{e}^{\Im(z)^2/4} U(a,\Re(z)).   \nonumber
\end{align}
Combining the above two results yields
\begin{align}
|\Phi_n(z)| &\leq \frac{2^{n/2} n! \mathrm{e}^{-(x^2-\rho^2)/2}}{\sqrt{2\pi}} \mathrm{e}^{x^2/2} U\left(n+\frac{1}{2},\sqrt{2}\rho\right) \nonumber \\
&=
\frac{2^{n/2} n! \mathrm{e}^{\rho^2/2}}{\sqrt{2\pi}} U\left(n+\frac{1}{2},\sqrt{2}\rho\right) \nonumber \\
&= |\Phi_n(\mathrm{i}\rho)|. \nonumber
\end{align}
which proves the desired result \eqref{eq:PhiMax}. This ends the proof.
\end{proof}

\begin{theorem}\label{thm:HermPolyProj}
If $f$ is analytic in the strip $\mathcal{S}_{\rho}$ for
some $\rho>0$ and $|f(z)|\leq\mathcal{K}|z|^{\sigma}$ for some
$\sigma\in\mathbb{R}$ as $|z|\rightarrow\infty$ within the strip, and if
\begin{equation}
V := \int_{\partial\mathcal{S}_{\rho}} |f(z)| |\mathrm{d}z| < \infty, \nonumber
\end{equation}
then, for $n\geq\max\{\lfloor\sigma\rfloor+1,0\}$, the following statements hold.
\begin{itemize}
\item[\rm(i)] 
The Hermite coefficient satisfies
\begin{equation}\label{eq:HermCoefBound}
|a_n| \leq \frac{V}{2^{n+1}\pi} U\left(\frac{n+1}{2},\frac{1}{2},\rho^2 \right).
\end{equation}

\item[\rm(ii)] The error of Hermite projection in the weighted $L^2$-norm satisfies
\begin{equation}\label{eq:HermProjBound}
\|f-\Pi_n^\mathrm{P}f\|_{L_{\omega}^2(\mathbb{R})} \leq \mathcal{K} \mathrm{e}^{-\rho\sqrt{2n}},
\end{equation}
where $\mathcal{K} \sim \mathrm{e}^{\rho^2/2}V/(2\sqrt{\pi\rho})$ for $n\gg1$. Hereafter, the symbol ``$A \sim B$'' for $n\gg1$ means $A=B(1+o(1))$ as $n\rightarrow\infty$.
\end{itemize}
\end{theorem}
\begin{proof}
By Theorem \ref{thm:HermCoeff}, Lemma \ref{lem:MaxPhi} and \eqref{eq:PhiFormula1}, we obtain immediately that
\begin{align}
|a_n| \leq \frac{|\Phi_n(\mathrm{i}\rho)|}{\gamma_n} \int_{\partial{\mathcal{S}}_{\rho}} |f(z)| |\mathrm{d}z|
= \frac{V}{2^{n+1}\pi} U\left(\frac{n+1}{2},\frac{1}{2},\rho^2\right),  \nonumber
\end{align}
which leads to \eqref{eq:HermCoefBound}. To show \eqref{eq:HermProjBound}, it is readily seen from \eqref{eq:HermCoefBound} that
\begin{align}
\|f-\Pi_n^\mathrm{P}f\|_{L_{\omega}^2(\mathbb{R})}^2 &=
\sum_{k=n+1}^{\infty} a_k^2 \gamma_k \leq \frac{V^2}{4\pi^{3/2}} \sum_{k=n+1}^{\infty} \frac{k!}{2^k} U^2\left(\frac{k+1}{2},\frac{1}{2},\rho^2\right). \nonumber
\end{align}
Furthermore, by \eqref{eq:AsyU} and the duplication formula and ratio asymptotics for the gamma function (see \cite[Equations (5.5.5) and (5.11.13)]{Olver2010}), it follows that, as $k\rightarrow\infty$,
\begin{align}
&\frac{k!}{2^k} U^2\left(\frac{k+1}{2},\frac{1}{2},\rho^2\right) \sim \frac{k!}{2^k} \frac{\pi \mathrm{e}^{\rho^2-2\rho\sqrt{2(k+1)}}}{\Gamma((k+1)/2) \Gamma((k+3)/2)}   \nonumber \\
&= \frac{\Gamma((k+2)/2)}{\Gamma((k+3)/2)} \sqrt{\pi} \mathrm{e}^{\rho^2-2\rho\sqrt{2(k+1)}} \sim \frac{2\sqrt{\pi}}{\sqrt{2(k+1)}} \mathrm{e}^{\rho^2-2\rho\sqrt{2(k+1)}}. \nonumber
\end{align}
Combining the above two estimates yields
\begin{align}
\|f-\Pi_n^\mathrm{P}f\|_{L_{\omega}^2(\mathbb{R})}^2 &\leq \mathcal{K} \sum_{k=n+1}^{\infty} \frac{\mathrm{e}^{-2\rho\sqrt{2(k+1)}}}{\sqrt{2(k+1)}} \nonumber
\end{align}
with $\mathcal{K} \sim \mathrm{e}^{\rho^2}V^2/(2\pi)$ for $n\gg1$. Since
\begin{equation}
\sum_{k=n+1}^{\infty} \frac{\mathrm{e}^{-2\rho\sqrt{2(k+1)}}}{\sqrt{2(k+1)}}
\leq \int_{n+1}^{\infty} \frac{\mathrm{e}^{-2\rho\sqrt{2x}}}{\sqrt{2x}}
\mathrm{d}x = \frac{\mathrm{e}^{-2\rho\sqrt{2(n+1)}}}{2\rho}, \nonumber
\end{equation}
the desired result \eqref{eq:HermProjBound} follows.
\end{proof}

\begin{remark}
The upper bound \eqref{eq:HermCoefBound} for Hermite coefficients is actually sharp. To see this, we consider the function $f(x)=1/(x^2+\tau^2)$ which is analytic in the strip $\mathcal{S}_{\rho}$ with $\rho=\tau-\epsilon$ and $\epsilon>0$ being arbitrarily close to zero. By Theorem \ref{thm:HermCoeff} and the residue theorem, one can show that the associated Hermite coefficients are explicitly given by
\begin{equation}
a_n = \frac{\cos(n\pi/2)}{2^n\tau} U\left(\frac{n+1}{2},\frac{1}{2}, \tau^2\right), \quad n\in\mathbb{N}_0. \nonumber
\end{equation}
Comparing this with \eqref{eq:HermCoefBound}, we see that the bound cannot be improved more than a constant factor.
\end{remark}

\begin{remark}
By \eqref{eq:AsyU} and the Stirling's formula, we have that, as $n\to \infty$,
\begin{align}
\frac{V}{2^{n+1}\pi} U\left(\frac{n+1}{2},\frac{1}{2},\rho^2 \right) &\sim \frac{\mathrm{e}^{\rho^2/2} V}{\sqrt{\pi}} \frac{\mathrm{e}^{-\rho\sqrt{2(n+1)}}}{2^{n+1} \Gamma((n+1)/2) \sqrt{(n+1)/2}}   \nonumber \\
&\sim \frac{\mathrm{e}^{\rho^2/2} V}{\sqrt{2}\pi} \frac{\mathrm{e}^{-\rho\sqrt{2(n+1)}}}{2^{n+1}} \left( \frac{2\mathrm{e}}{n+1} \right)^{(n+1)/2}. \nonumber
\end{align}
We then conclude from \eqref{eq:HermCoefBound} that the Hermite coefficient $a_n$ decays at least at a super-exponential rate.
\end{remark}

\begin{remark}
Convergence of Hermite polynomial approximation for analytic functions in the weighted norm has been studied in different contexts, e.g., stochastic collocation method \cite{Babuska2007} and deep neural network \cite{Schwab2023}. For example, it was recently shown in \cite[Theorem 4.5]{Schwab2023} that the Hermite projection of degree $n$ using the scaled Hermite polynomials $\{H_k(x/\sqrt{2})\}_{k=0}^{\infty}$ converges at the rate $O(\mathrm{e}^{-\rho\sqrt{n/2}})$. In contrast, our analysis gives the convergence rate $O(\mathrm{e}^{-\rho\sqrt{n}})$ after extending \eqref{eq:HermProjBound} to this case, which is clearly better.
\end{remark}

\subsection{Projections using Hermite functions}
For $n\in\mathbb{N}_0$, the Hermite functions are defined by
\begin{equation}\label{def:HermFunc}
\psi_n(x) = \frac{H_n(x)}{\sqrt{\gamma_n}} \mathrm{e}^{-x^2/2} ,
\end{equation}
which are preferable in physical problems, since they are the eigenfunctions of the quantum harmonic oscillator. It is known that the Hermite functions $\{\psi_n\}_{n=0}^{\infty}$ form a complete orthonormal basis for the Hilbert space $L^2(\mathbb{R})$ and by \eqref{eq:HermBound},
\begin{equation}\label{eq:PsiBound}
|\psi_n(x)| \leq \pi^{-1/4}, \qquad  \forall x\in\mathbb{R}.
\end{equation}
We now introduce the space $\mathbb{H}_n=\mathrm{span}\{\psi_k\}_{k=0}^{n}$, which is spanned
by Hermite functions, and let
$\Pi_n^{\mathrm{F}}:L^2(\mathbb{R})\rightarrow\mathbb{H}_n$ be the
orthogonal projection operator from $L^2(\mathbb{R})$ upon
$\mathbb{H}_n$. For any $f\in{L}^2(\mathbb{R})$, its orthogonal projection upon
$\mathbb{H}_n$ reads
\begin{equation}\label{eq:HermFuncExpan}
(\Pi_n^{\mathrm{F}}f)(x) = \sum_{k=0}^{n} c_k \psi_k(x), \qquad  c_k
= \int_{\mathbb{R}}f(x)\psi_k(x) \mathrm{d}x.
\end{equation}
Let $\|\cdot\|_{L^2(\mathbb{R})}$ and $\|\cdot\|_{L^{\infty}(\mathbb{R})}$ stand for the norms induced by the usual inner product and by the maximum norm over $\mathbb{R}$, respectively. Below we give some error bounds with asymptotically explicit constants for $\Pi_n^{\mathrm{F}}$ both in the $L^2$-norm and in the $L^{\infty}$-norm.
\begin{theorem}\label{thm:ProHermFunc}
If $f$ is analytic in the strip $\mathcal{S}_{\rho}$ for
some $\rho>0$ and $|\mathrm{e}^{z^2/2}f(z)|\leq\mathcal{K}|z|^{\sigma}$ for some
$\sigma\in\mathbb{R}$ as $|z|\rightarrow\infty$ within the strip, and if
\begin{equation}
\widehat{V} := \int_{\partial\mathcal{S}_{\rho}} |\mathrm{e}^{z^2/2} f(z)| |\mathrm{d}z| < \infty, \nonumber
\end{equation}
then, for $n\geq\max\{\lfloor\sigma\rfloor+1,0\}$, the following statements hold.
\begin{itemize}
\item[\rm(i)] 
The Hermite coefficient in \eqref{eq:HermFuncExpan} satisfies
\begin{equation}\label{eq:Herm3}
|c_n| \leq \frac{\widehat{V}\sqrt{\gamma_n}}{ 2^{n+1}\pi} U\left(\frac{n+1}{2},\frac{1}{2},\rho^2 \right),
\end{equation}
and, for $n\gg1$, a more explicit bound is given by
\begin{equation}\label{eq:estCn}
|c_n| \leq \mathcal{K} n^{-1/4} \mathrm{e}^{-\rho\sqrt{2n}},
\end{equation}
where $\mathcal{K} \sim \mathrm{e}^{\rho^2/2}\widehat{V}/(2^{3/4}\sqrt{\pi})$.

\item[\rm(ii)] The error of Hermite projection in the $L^2$-norm satisfies
\begin{equation}\label{eq:Herm4}
\|f-\Pi_n^\mathrm{F}f\|_{L^2(\mathbb{R})} \leq
\mathcal{K} \mathrm{e}^{-\rho\sqrt{2n}},
\end{equation}
where $\mathcal{K}\sim
\mathrm{e}^{\rho^2/2}\widehat{V}/(2\sqrt{\pi\rho})$ for $n\gg1$.

\item[\rm(iii)] The error of Hermite projection in the $L^{\infty}$-norm satisfies
\begin{equation}\label{eq:Herm5}
\|f-\Pi_n^\mathrm{F}f\|_{L^{\infty}(\mathbb{R})} \leq
\mathcal{K} n^{1/4} \mathrm{e}^{-\rho\sqrt{2n}},
\end{equation}
where $\mathcal{K} \sim \mathrm{e}^{\rho^2/2}\widehat{V}/(2^{1/4}\pi^{3/4}\rho)$ for $n\gg1$.
\end{itemize}
\end{theorem}
\begin{proof}
By \eqref{def:HermFunc} and \eqref{eq:HermFuncExpan}, we have
\[
\frac{c_n}{\sqrt{\gamma_n}} = \frac{1}{\gamma_n} \int_{\mathbb{R}} f(x) H_n(x) \mathrm{e}^{-x^2/2} \mathrm{d}x.
\]
Note that the right hand side above is exactly the $n$-th expansion coefficient of $\mathrm{e}^{x^2/2}f(x)$ in terms of Hermite polynomials, the estimate \eqref{eq:Herm3} then follows directly from \eqref{eq:HermCoefBound}. Furthermore, \eqref{eq:estCn} follows from a combination of \eqref{eq:Herm3} and \eqref{eq:AsyU}. As for \eqref{eq:Herm4}, by \eqref{eq:estCn} and the orthogonality of $\{\psi_n\}$, it is readily seen that
\begin{align}
\|f-\Pi_n^\mathrm{F}f\|_{L^2(\mathbb{R})}^2 &= \sum_{k=n+1}^{\infty}
c_k^2 \leq \mathcal{K} \sum_{k=n+1}^{\infty} \frac{\mathrm{e}^{-2\rho\sqrt{2k}}}{\sqrt{k}}, \nonumber
\end{align}
where $\mathcal{K}\sim\mathrm{e}^{\rho^2}\widehat{V}^2/(2^{3/2}\pi)$ for $n\gg1$. Note that
\begin{equation}
\sum_{k=n+1}^{\infty} \frac{\mathrm{e}^{-2\rho\sqrt{2k}}}{\sqrt{k}}
\leq \int_{n}^{\infty} \frac{\mathrm{e}^{-2\rho\sqrt{2x}}}{\sqrt{x}}
\mathrm{d}x = \frac{\mathrm{e}^{-2\rho\sqrt{2n}}}{\sqrt{2}\rho}. \nonumber
\end{equation}
The error bound \eqref{eq:Herm4} follows by combining the above two inequalities.

To show \eqref{eq:Herm5}, we see from the inequalities \eqref{eq:PsiBound} and \eqref{eq:estCn} that
\begin{align}
\|f-\Pi_n^\mathrm{F}f\|_{L^{\infty}(\mathbb{R})} &\leq
\frac{1}{\pi^{1/4}} \sum_{k=n+1}^{\infty} |c_k| \leq \mathcal{K} \sum_{k=n+1}^{\infty}
\frac{\mathrm{e}^{-\rho\sqrt{2k}}}{k^{1/4}}, \nonumber 
\end{align}
where $\mathcal{K} \sim  \mathrm{e}^{\rho^2/2}\widehat{V}/(2\pi)^{3/4}$ for $n\gg1$. For the summation in the last inequality, we note that
\begin{equation}
\sum_{k=n+1}^{\infty} \frac{\mathrm{e}^{-\rho\sqrt{2k}}}{k^{1/4}}
\leq \int_{n}^{\infty} \frac{\mathrm{e}^{-\rho\sqrt{2x}}}{x^{1/4}}
\mathrm{d}x = \frac{2^{1/4}}{\rho^{3/2}}
\Gamma\left(\frac{3}{2},\rho\sqrt{2n}\right), \nonumber
\end{equation}
where $\Gamma(a,z)$ is the incomplete gamma function. The bound \eqref{eq:Herm5} follows by combining
the fact that $\Gamma(a,z)=z^{a-1}\mathrm{e}^{-z}(1+O(z^{-1}))$ as $z\to\infty$ (see \cite[Equation~(8.11.2)]{Olver2010})
with the above two inequalities. This ends the proof of Theorem \ref{thm:ProHermFunc}.
\end{proof}

\begin{remark}
Hille in \cite{Hille1939,Hille1940} showed that $|c_n| = O(\mathrm{e}^{-\rho\sqrt{2n}})$ if $f$ is analytic in the strip $\mathcal{S}_{\rho}$ and $|f(x+\mathrm{i}y)|\leq \mathcal{K}\exp(-|x|\sqrt{\rho^2-y^2})$ for $x+\mathrm{i}y\in\mathcal{S}_{\rho}$. Boyd in \cite{Boyd1984} analyzed the function $f(x)=\mathrm{e}^{-x^2}/(x^2+\rho^2)$ and showed that its Hermite coefficients satisfy $|c_n|\sim 2^{5/4}\sqrt{\pi}R n^{-1/4} \mathrm{e}^{-\rho\sqrt{2n}}$ as $n\rightarrow\infty$, where $R$ is the residue of $f$ at the poles $z=\pm\mathrm{i}\rho$. Our bound \eqref{eq:estCn} contains the prefactor $n^{-1/4}$ and is more precise than Hille's result, but with a stronger assumption on the asymptotic behavior of the underlying function at infinity. Moreover, our result is more general than Boyd's result in the sense that it has no restriction on the types of singularities.
\end{remark}

\begin{remark}
Let $\mathcal{S}$ be the Schwartz space of rapidly decaying $C^{\infty}$ functions on $\mathbb{R}$. For every integer $s\leq n+1$ and every $f\in\mathcal{S}$, it has been proved in \cite[Chapter III]{Lubich2008} that
\[
\|f-\Pi_n^\mathrm{F}f\|_{L^2(\mathbb{R})} \leq \frac{1}{\sqrt{(n+1)n\cdots(n-s+2)}} \|A^sf\|_{L^2(\mathbb{R})},
\]
where $A$ is the Dirac's ladder operator defined by $Af = 2^{-1/2}(x + \mathrm{d}/\mathrm{d}x)f$. When $f$ is analytic in the strip $\mathcal{S}_{\rho}$, it is clear that this result gives a superalgebraic rate of convergence of $\Pi_n^\mathrm{F}f$ in the $L^2$-norm, which is inferior to our result \eqref{eq:Herm4}.
\end{remark}

\begin{figure}[htbp]
\centering
\includegraphics[width=0.49\textwidth]{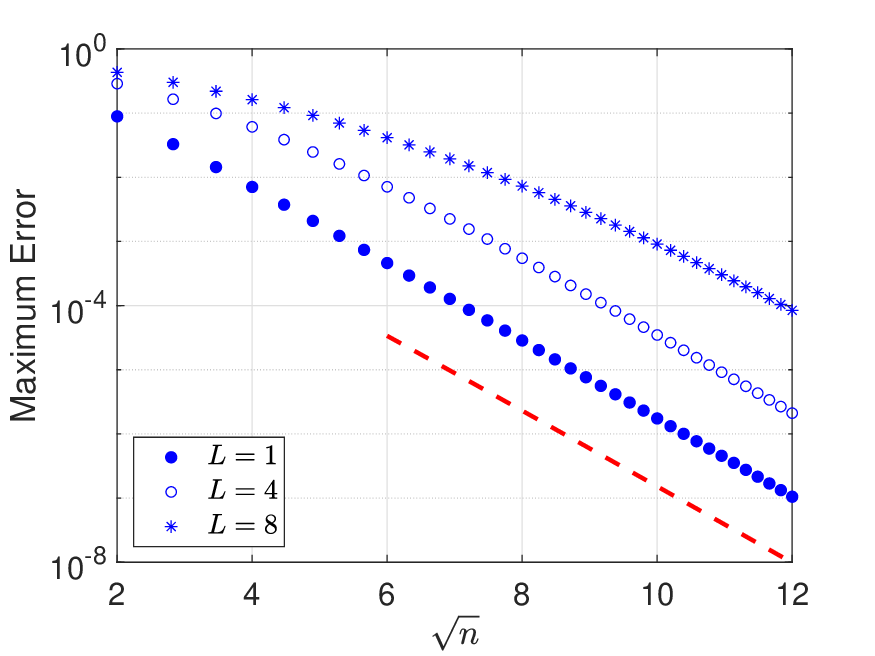}~
\includegraphics[width=0.49\textwidth]{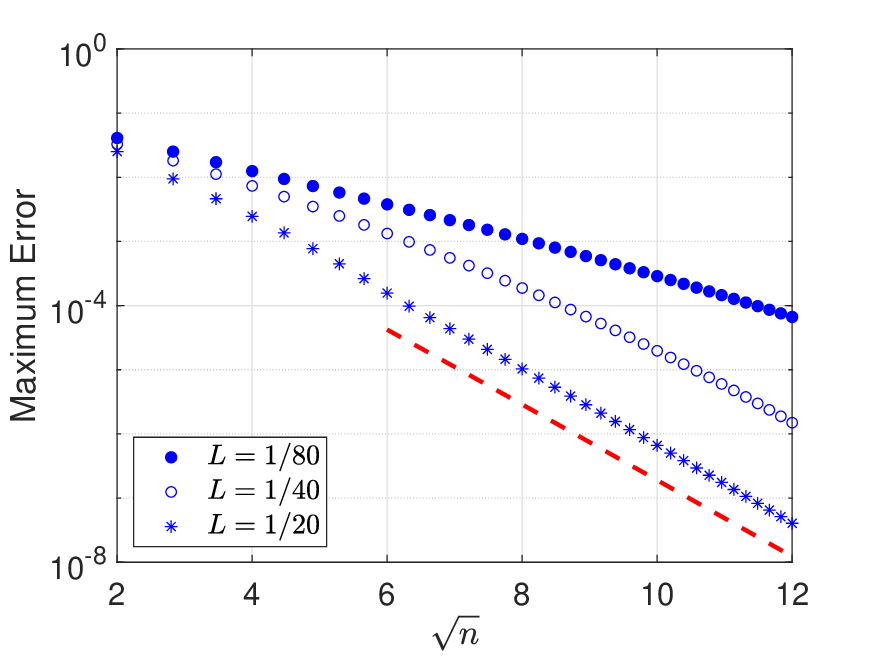}
\caption{Maximum error of $\Pi_n^\mathrm{F}f$ as a function of $\sqrt{n}$ for $f(x) = \mathrm{e}^{-Lx^2}/(x^2+1)$. The dashed lines show the predicted rate $O(n^{1/4}\mathrm{e}^{-\sqrt{2n}})$.}\label{fig:MaxErr2}
\end{figure}

We present numerical illustrations of Theorem \ref{thm:ProHermFunc}. Consider the function
$f(x) = \mathrm{e}^{-Lx^2}/(x^2+1)$ and $L>0$. Clearly, this function is analytic in the strip $\mathcal{S}_{\rho}$ with $\rho=1-\epsilon$, where $\epsilon>0$ is arbitrarily close to zero. From \eqref{eq:Herm5} we conclude that the convergence rate of $\Pi_n^\mathrm{F}f$ in the $L^{\infty}$-norm is $O(n^{1/4}\mathrm{e}^{-\sqrt{2n}})$ when $L\geq1/2$. In Figure \ref{fig:MaxErr2} we plot the maximum error of $\Pi_n^\mathrm{F}f$ for several different values of $L$. We see that the convergence rate of $\Pi_n^\mathrm{F}f$ matches the predicted rate quite well when $L$ is close to one. As $L\rightarrow\infty$ or $L\rightarrow0$, however, the accuracy of $\Pi_n^\mathrm{F}f$ deteriorates gradually. Indeed, a straightforward calculation shows that $\widehat{V}=O(\exp(L\rho^2))$ as $L\rightarrow\infty$ and therefore the error bound in \eqref{eq:Herm5} involves an exponentially large constant. For $L<1/2$, a straightforward calculation shows that $\widehat{V}=\infty$ and the error bound in \eqref{eq:Herm5} fails in this case. However, we see from the right panel of Figure \ref{fig:MaxErr2} that $\Pi_n^\mathrm{F}f$ for $L=1/20$ still converges at the rate $O(n^{1/4}\mathrm{e}^{-\sqrt{2n}})$, this implies that the conditions of Theorem \ref{thm:ProHermFunc} might be relaxed further.

\begin{figure}[htbp]
\centering
\includegraphics[width=0.49\textwidth]{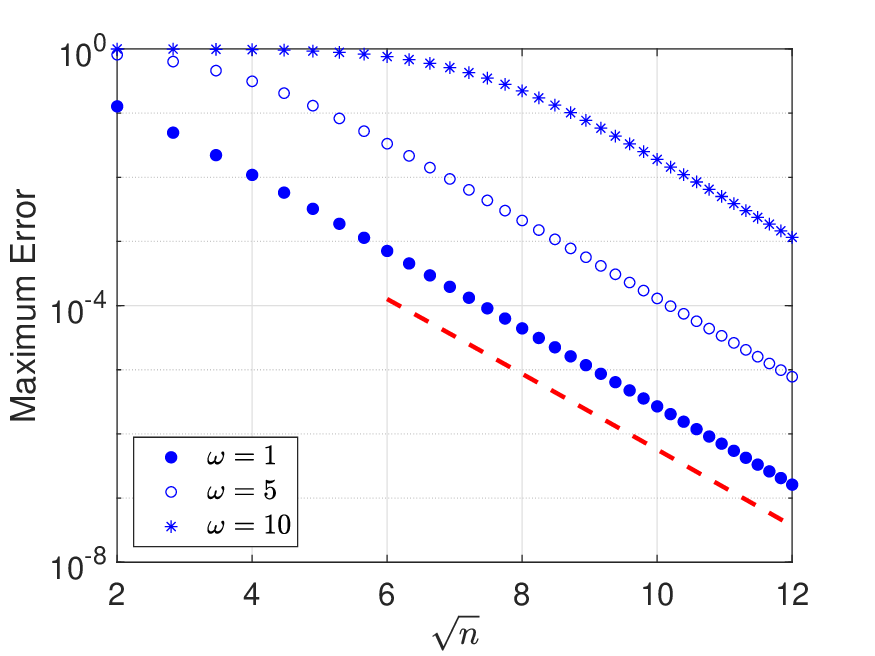}
\includegraphics[width=0.49\textwidth]{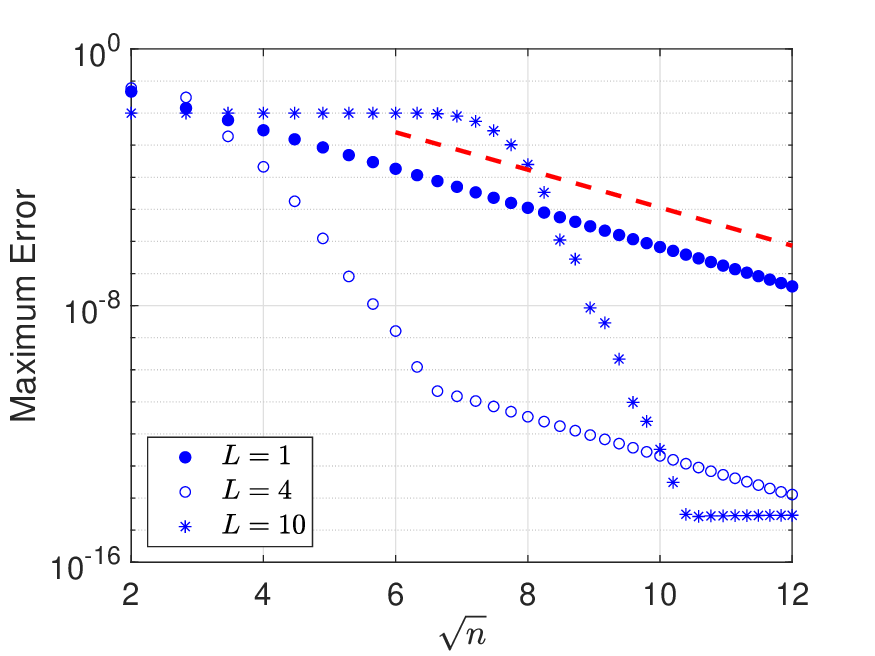}
\caption{Maximum error of $\Pi_n^\mathrm{F}f$ as a function of $\sqrt{n}$ for $f(x)=\mathrm{e}^{-x^2}\cos(\omega x)/(x^2+1)$ (left) and $f(x)=\mathrm{e}^{-(x-L)^2}/(x^2+1)$ (right). Here the dashed lines show the predicted rate $O(n^{1/4}\mathrm{e}^{-\sqrt{2n}})$.}\label{fig:MaxErr3}
\end{figure}

In Figure \ref{fig:MaxErr3} we plot the maximum error of $\Pi_n^\mathrm{F}f$ for $f(x)=\mathrm{e}^{-x^2}\cos(\omega x)/(x^2+1)$ and $f(x)=\mathrm{e}^{-(x-L)^2}/(x^2+1)$. In the former case, we see that the convergence rate of $\Pi_n^\mathrm{F}f$ for each fixed $\omega$ matches the predicted rate as $n\rightarrow\infty$. For fixed $n$, however, the accuracy of $\Pi_n^\mathrm{F}f$ deteriorates gradually as $\omega$ increases. This is due to the fact that $\widehat{V}=O(\exp(\omega\rho))$ as $\omega\rightarrow\infty$ and hence the error bound in \eqref{eq:Herm5} involves an exponentially large constant. For the latter case, we see that the convergence rate of $\Pi_n^\mathrm{F}f$ matches the predicted rate for almost all $n$ when $L=1$ and for $n\geq44$ when $L=4$. When $L=10$, however, the error of $\Pi_n^\mathrm{F}f$ retains at the level of $10^{-2}$ for $n\leq48$ and then decays at a much faster rate until it reaches the level of $10^{-15}$. For $n\leq48$, this might follow from the fact that $\widehat{V}=O(\exp(L^2))$ as $L\rightarrow\infty$, but we still have no idea to explain the faster convergence rate for $n\geq48$.

\section{Convergence analysis of Hermite spectral interpolations}\label{sec:Interp}
In this section, we are concerned with interpolation methods using Hermite polynomials and Hermite functions, which are also known as the Hermite spectral interpolation\footnote{Note that Hermite interpolation is often referred as polynomial interpolation
with derivative conditions.}. Let $\{x_j\}_{j=1}^{n}$ be the
zeros of $H_{n}(x)$ and we assume that they are arranged in
ascending order, i.e., $-\infty<x_1<\cdots<x_n<\infty$. By the symmetry relation \eqref{eq:HermSymm} of
Hermite polynomials, it is easily seen that $x_j=-x_{n-j+1}$ for $j=1,\ldots,n$. Let $p_{n}\in\mathbb{P}_{n-1}$ be the unique polynomial which interpolates $f(x)$ at the points
$\{x_j\}_{j=1}^{n}$, i.e.,
\begin{equation}\label{eq:HermiteInter}
p_{n}(x_j) = f(x_j), \quad  j=1,\ldots,n,
\end{equation}
and let $h_{n}\in\mathbb{H}_{n-1}$ be the unique function which interpolates $f(x)$ at the points
$\{x_j\}_{j=1}^{n}$, i.e.,
\begin{equation}\label{eq:hn}
h_{n}(x_j) = f(x_j), \qquad j=1,\ldots,n.
\end{equation}
The convergence analysis of $p_{n}(x)$ and $h_n(x)$ will be the main topic of this section.

We start with presenting a contour integral representation of the remainder for Hermite spectral interpolation by extending the idea of \cite[Lemma~4.1]{Wang2023c}.
\begin{lemma}\label{lem:ContourInterp}
If $f$ is analytic in the strip $\mathcal{S}_{\rho}$ for
some $\rho>0$ and $|f(z)|\leq\mathcal{K}|z|^{\sigma}$ for some
$\sigma\in\mathbb{R}$ as $|z|\rightarrow\infty$ within the strip, then for any $x\in\mathbb{R}$ and $n\geq\max\{\lfloor\sigma\rfloor+1,0\}$, we have
\begin{equation}
f(x) - p_{n}(x) =
\frac{1}{2\pi\mathrm{i}}\int_{\partial{\mathcal{S}}_{\rho}}
\frac{H_{n}(x)f(z)}{H_{n}(z)(z-x)} \mathrm{d}z,
\end{equation}
where $p_{n}$ is the polynomial of degree $n-1$ determined through \eqref{eq:HermiteInter}.
\end{lemma}
\begin{proof}
Let $\eta>0$ be large enough such that $x\in(-\eta,\eta)$ and set
\begin{equation}
V_{\pm}:~~z=\pm\eta+\mathrm{i}y,\qquad y\in[-\rho,\rho]. \nonumber
\end{equation}
Note that $V_{\pm}$ are related to the contours $\Gamma_{\pm}^{\mathrm{R}}$ and $\Gamma_{\pm}^{\mathrm{L}}$ defined in \eqref{def:contours} by
$V_{+}=\Gamma_{+}^{\mathrm{R}}\cup\Gamma_{-}^{\mathrm{R}}$ and
$V_{-}=\Gamma_{+}^{\mathrm{L}}\cup\Gamma_{-}^{\mathrm{L}}$ with
$\varepsilon=0$. By Hermite's contour integral \cite[Theorem~3.6.1]{Davis1975}, the
remainder of $p_{n}$ can be written as
\begin{equation}\label{eq:HermRemainder}
f(x) - p_{n}(x) =
\frac{1}{2\pi\mathrm{i}} \int_{\Gamma_{-}^{\mathrm{U}}\cup\Gamma_{+}^{\mathrm{U}}\cup V_{+}\cup V_{-}}
\frac{H_{n}(x)f(z)}{H_{n}(z)(z-x)} \mathrm{d}z.
\end{equation}
For $z\in{V}_{+}$, we have
\begin{equation}
\left|\int_{V_{+}} \frac{H_{n}(x)f(z)}{H_{n}(z)(z-x)}
\mathrm{d}z\right| \leq
\left|\frac{H_{n}(x)}{H_{n}(\eta)}\right| \int_{V_{+}}
\frac{|f(z)|}{|z-x|} |\mathrm{d}z| \leq 2\rho \cdot \left|\frac{H_{n}(x)}{H_{n}(\eta)}\right| \max_{z\in{V_{+}}}\left|\frac{f(z)}{z-x}\right| ,  \nonumber
\end{equation}
and similarly,
\begin{equation}
\left|\int_{V_{-}} \frac{H_{n}(x)f(z)}{H_{n}(z)(z-x)}
\mathrm{d}z\right| \leq
\left|\frac{H_{n}(x)}{H_{n}(\eta)}\right| \int_{V_{-}}
\frac{|f(z)|}{|z-x|} |\mathrm{d}z| \leq 2\rho \cdot \left|\frac{H_{n}(x)}{H_{n}(\eta)}\right| \max_{z\in{V_{-}}}\left|\frac{f(z)}{z-x}\right|. \nonumber
\end{equation}
Recall that $|f(z)|\leq \mathcal{K}|z|^{\sigma}$ as $z\rightarrow\infty$ within the strip $\mathcal{S}_{\rho}$, we see that both the above two bounds behave like $O(\eta^{\sigma-n-1})$ as $\eta\rightarrow\infty$.
Since $n\geq\max\{\lfloor\sigma\rfloor+1,0\}$, it is easily seen that
\begin{equation*}
\sigma-n \left\{
           \begin{array}{ll}
             <0,           & \hbox{$\sigma>-1$,} \\[1ex]
             \leq -1,      & \hbox{$\sigma=-1$,} \\[1ex]
             \leq \sigma,  & \hbox{$\sigma<-1$.}
           \end{array}
         \right.
\end{equation*}
On the other hand, note that $f(z)/(H_n(z)(z-x))$ behaves like $O(z^{\sigma-n-1})$ as $z\rightarrow\infty$ within the strip $\mathcal{S}_{\rho}$, we deduce that the integrals
\[
\int_{\Gamma_{-}^{\mathrm{U}}} \frac{H_{n}(x)f(z)}{H_{n}(z)(z-x)} \mathrm{d}z, \quad  \int_{\Gamma_{+}^{\mathrm{U}}} \frac{H_{n}(x)f(z)}{H_{n}(z)(z-x)} \mathrm{d}z,
\]
converge as $\eta\rightarrow\infty$. Thus, by taking $\eta\rightarrow\infty$ on both sides of \eqref{eq:HermRemainder} and noting that
$\partial{\mathcal{S}}_{\rho}= \Gamma_{-}^{\mathrm{U}}\cup \Gamma_{+}^{\mathrm{U}}$, we obtain the desired result from the above three results. This ends the proof. 
\end{proof}

We next prove an interesting property of Hermite polynomials, which might be of independent interest.
\begin{lemma}\label{lem:MinHerm}
Let $x,y\in\mathbb{R}$ with $y\neq0$. For each $n\in\mathbb{N}$, it holds that
\begin{equation}\label{eq:minHerm}
\min_{x\in\mathbb{R}}|H_{n}(x+\mathrm{i}y)| \geq |H_{n}(\mathrm{i}y)| \left\{
           \begin{array}{ll}
             1,              & \hbox{$n$ odd,} \\[2ex]
           {\displaystyle  \sqrt{1 - \left|\frac{H_n(0)}{H_{n}(\mathrm{i}y)}\right|^2}}, & \hbox{$n$ even,}
           \end{array}
         \right.
\end{equation}
and for any $\rho>0$,
\begin{align}\label{eq:minHerm2}
\min_{z\in\partial\mathcal{S}_{\rho}}|H_n(z)| &\geq |H_{n}(\mathrm{i}\rho)| \left\{
           \begin{array}{ll}
             1,              & \hbox{$n$ odd,} \\[2ex]
           {\displaystyle  \sqrt{1 - \left|\frac{H_n(0)}{H_{n}(\mathrm{i}\rho)}\right|^2}}, & \hbox{$n$ even,}
           \end{array}
         \right.
\end{align}
and
\begin{align}\label{eq:minHerm3}
\min_{z\in\partial\mathcal{S}_{\rho}}|H_n(z)| &= |H_n(\mathrm{i}\rho)| \left\{
           \begin{array}{ll}
             1,              & \hbox{$n$ odd,} \\[2ex]
           {\displaystyle  1 + O\left(\mathrm{e}^{-2\rho\sqrt{2n}}\right)}, & \hbox{$n$ even,}
           \end{array}
         \right.
\end{align}
where the case for even $n$ in the last equation is understood for large $n$.
\end{lemma}
\begin{proof}
We first consider \eqref{eq:minHerm}. Let $P\in\mathbb{P}_n$ be a polynomial with $n$ real roots, by \cite[Theorem~2.1]{Patrick1971} we know that
\begin{equation}\label{eq:Pequality}
   |P(x+\mathrm{i}y)|^2 = \sum_{k=0}^{n} L_k(P;x) y^{2k},
\end{equation}
where $x,y\in\mathbb{R}$ and
\[
L_k(P;x) = \sum_{j=0}^{2k} \frac{(-1)^{j+k}}{(2k)!} \binom{2k}{j} P^{(j)}(x) P^{(2k-j)}(x).
\]
Moreover, by \cite[Theorem~2.2]{Patrick1971} we also know that $L_k(P;x)\geq0$ for all $x\in\mathbb{R}$. Below we consider the case of Hermite polynomials, i.e., $P(x)=H_n(x)$. In this case, it was proved in \cite[Theorem~2]{Alex2012} that each $L_k(H_n;x)$, $k=1,\ldots,n$, is monotonously decreasing on $(-\infty,0]$ and monotonously increasing on $[0,\infty)$. Therefore, for $x\in\mathbb{R}$, we have from \eqref{eq:Pequality} that
\begin{align}
|H_n(x+\mathrm{i}y)|^2 &= (H_n(x))^2 + \sum_{k=1}^{n} L_k(H_n;x) y^{2k} \nonumber \\
&\geq (H_n(x))^2 + \sum_{k=1}^{n} L_k(H_n;0) y^{2k} \nonumber \\
&= (H_n(x))^2 + |H_n(\mathrm{i}y)|^2 - (H_n(0))^2 \nonumber \\
&\geq |H_n(\mathrm{i}y)|^2 - (H_n(0))^2, \nonumber
\end{align}
where we have used the fact that $L_0(P;x)=(P(x))^2$ in the first step. The desired result \eqref{eq:minHerm} for odd $n$ follows immediately from the above equation and the fact that $H_n(0)=0$. As for even $n$, since the leading coefficient of $H_n(z)$ is $2^n$, we can write $H_n(z)=2^n(z-x_1)\cdots(z-x_n)$ where $\{x_j\}_{j=1}^{n}$ are the zeros of $H_n(z)$, and thus
\[
|H_n(\mathrm{i}y)| = 2^n \sqrt{(y^2+x_1^2)\cdots(y^2+x_n^2)}. 
\]
Clearly, $|H_n(0)|<|H_n(\mathrm{i}y)|$ for $y\neq0$, hence the desired result \eqref{eq:minHerm} for even $n$ follows. As for \eqref{eq:minHerm2}, it follows from \eqref{eq:minHerm} and the fact that $|H_n(-\mathrm{i}\rho)|=|H_n(\mathrm{i}\rho)|$. As for \eqref{eq:minHerm3}, it remains to show the case of even $n$. By \eqref{eq:HermFunAsy} and the fact that $|H_n(0)|=\Gamma(n+1)/\Gamma(n/2+1)$, it is readily seen that
\[
\left|\frac{H_n(0)}{H_n(\mathrm{i}\rho)} \right| = 2\mathrm{e}^{\rho^2/2} \mathrm{e}^{-\rho\sqrt{2n+1}} \left(1 + O(n^{-1/2})\right), \quad n\rightarrow\infty.
\]
Thus, the case of even $n$ in \eqref{eq:minHerm3} follows by combining the above result with \eqref{eq:minHerm2}. This ends the proof.
\end{proof}

\begin{figure}[htbp]
\centering
\includegraphics[width=0.49\textwidth,height=0.45\textwidth]{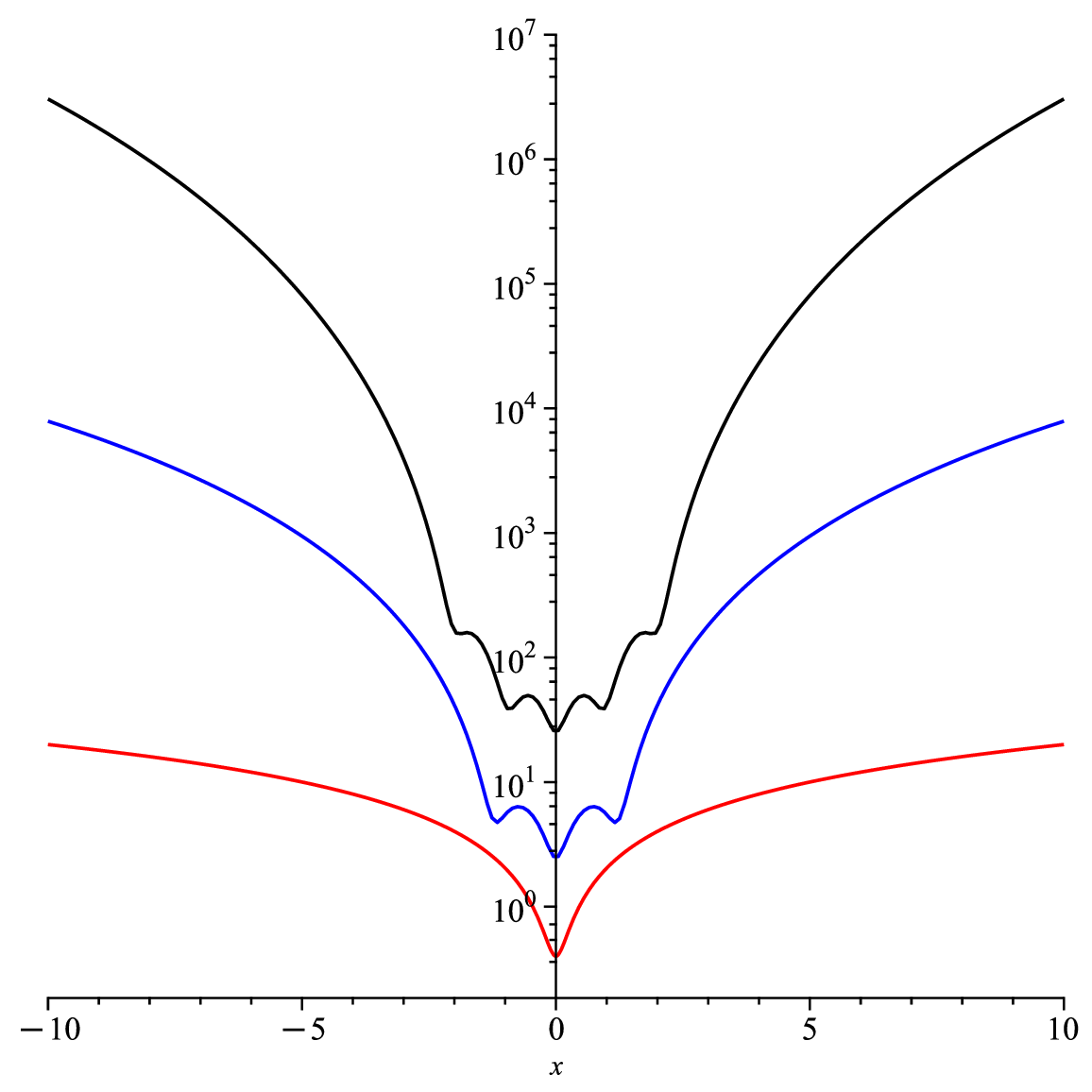}~
\includegraphics[width=0.49\textwidth,height=0.45\textwidth]{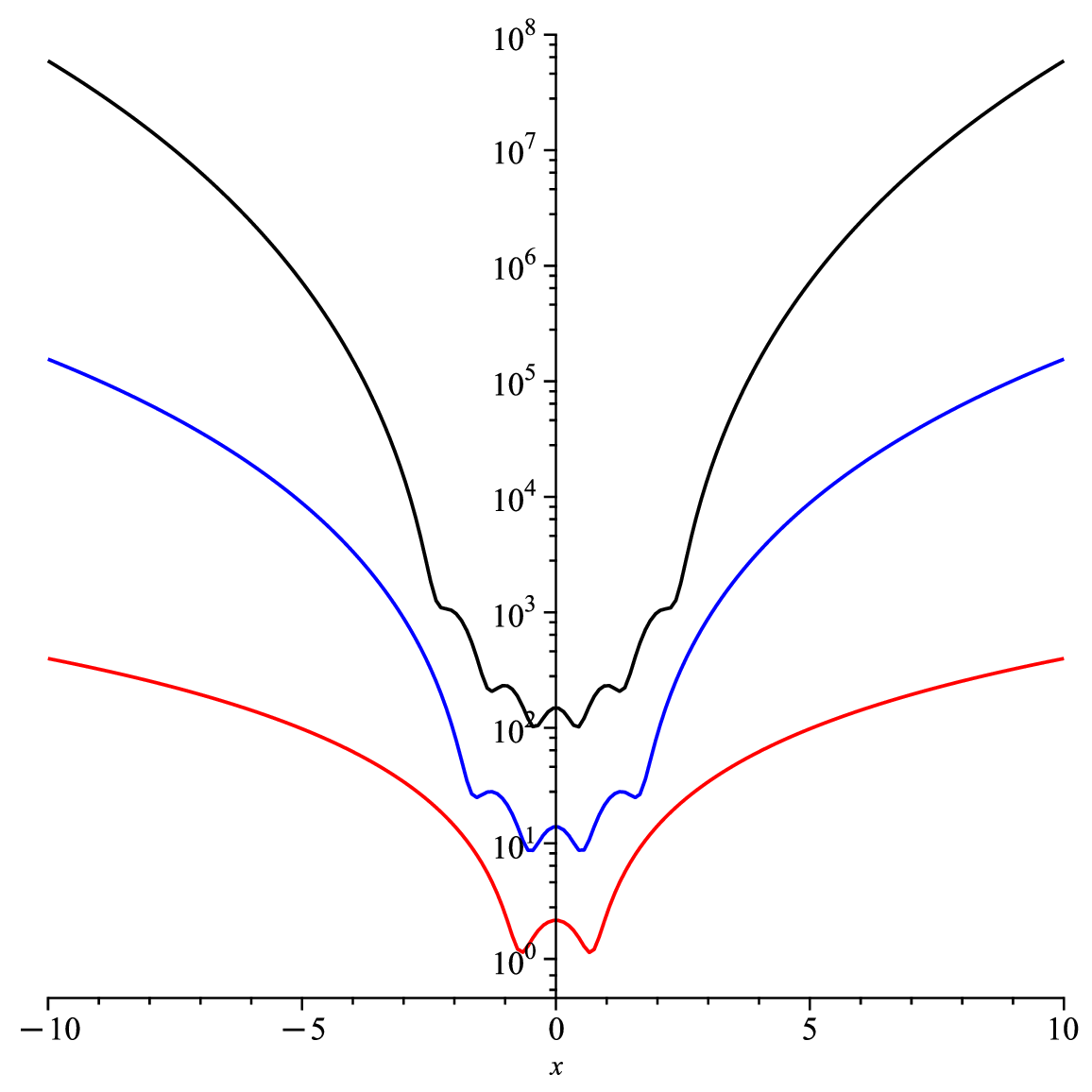}
\caption{Plot of $|H_n(x+\mathrm{i}\rho)|$ for $\rho=1/5$ and $x\in[-10,10]$. Here $n=1,3,5$ (left) and $n=2,4,6$ (right).}\label{fig:LocMin}
\end{figure}

\begin{remark}
We see that $|H_n(z)|$ for $z\in\partial\mathcal{S}_{\rho}$ attains its minimum at $z=\pm\mathrm{i}\rho$ for odd $n$. When $n$ is even, however, the locations where $|H_n(z)|$ attains its minimum for $z\in\partial\mathcal{S}_{\rho}$ will slightly deviate from $\pm\mathrm{i}\rho$ since $H_n(0)\neq0$ and those locations will approach to $\pm\mathrm{i}\rho$ as $n$ increases; see Figure \ref{fig:LocMin} for an illustration.
\end{remark}

\begin{theorem}
If $f$ is analytic in the strip $\mathcal{S}_{\rho}$ for
some $\rho>0$ and $|f(z)|\leq\mathcal{K}|z|^{\sigma}$ for some
$\sigma\in\mathbb{R}$ as $|z|\rightarrow\infty$ within the strip and if $V<\infty$, then for $n\geq\max\{\lfloor\sigma\rfloor+1,0\}$, we have
\begin{equation}\label{eq:errorHSI0}
\|f - p_{n} \|_{L_{\omega}^2(\mathbb{R})} \leq  \frac{V}{2\pi\rho} \left\{
           \begin{array}{ll}
           {\displaystyle \frac{\sqrt{\gamma_n}}{|H_n(\mathrm{i}\rho)|}},              & \hbox{$n$ odd,} \\[3ex]
           {\displaystyle \frac{\sqrt{\gamma_n}}{\sqrt{|H_n(\mathrm{i}\rho)|^2-|H_n(0)|^2}}}, & \hbox{$n$ even.}
           \end{array}
         \right.
\end{equation}
Moreover, for $n\gg1$, we have the following more explicit error bound
\begin{equation}\label{eq:errorHSI}
\|f - p_{n} \|_{L_{\omega}^2(\mathbb{R})} \leq \mathcal{K} n^{1/4}
\mathrm{e}^{-\rho\sqrt{2n}},
\end{equation}
where $\mathcal{K} \sim \mathrm{e}^{\rho^2/2}V/(2^{1/4}\sqrt{\pi}\rho)$.
\end{theorem}
\begin{proof}
By Lemma \ref{lem:ContourInterp}, we have
\begin{align}
\|f - p_{n} \|_{L_{\omega}^2(\mathbb{R})} &= \sqrt{\int_{\mathbb{R}}
\omega(x) |f(x) - p_{n}(x)|^2 \mathrm{d}x} \leq \frac{\sqrt{\gamma_{n}}}{2\pi\rho}
\int_{\partial{\mathcal{S}}_{\rho}} \frac{|f(z)|}{|H_{n}(z)|}
|\mathrm{d}z|,  \nonumber
\end{align}
where we have used the fact $\min|z-x|=\rho$ for $z\in\partial{\mathcal{S}}_{\rho}$ and $x\in\mathbb{R}$. The error bound \eqref{eq:errorHSI0} then follows immediately from the last inequality and \eqref{eq:minHerm2}. As for \eqref{eq:errorHSI}, it follows directly by combining \eqref{eq:errorHSI0} with \eqref{eq:HermFunAsy} and the duplication formula and ratio asymptotics for the gamma functions (see \cite[Equations (5.5.5) and (5.11.13)]{Olver2010}). We omit the details here.
\end{proof}

Finally, we consider the convergence of $h_{n}(x)$ in the $L^{\infty}$-norm.
\begin{theorem}
If $f$ is analytic in the strip $\mathcal{S}_{\rho}$ for
some $\rho>0$ and $|\mathrm{e}^{z^2/2}f(z)|\leq\mathcal{K}|z|^{\sigma}$ for some
$\sigma\in\mathbb{R}$ as $|z|\rightarrow\infty$ within the strip and if $\widehat{V}<\infty$, then for $n\geq\max\{\lfloor\sigma\rfloor+1,0\}$, we have
\begin{equation}\label{eq:errorHSI1}
\|f - h_{n} \|_{L^{\infty}(\mathbb{R})} \leq \frac{\widehat{V}}{2\pi^{5/4}\rho } \left\{
           \begin{array}{ll}
           {\displaystyle \frac{\sqrt{\gamma_n}}{|H_n(\mathrm{i}\rho)|}},              & \hbox{$n$ odd,} \\[3ex]
           {\displaystyle \frac{\sqrt{\gamma_n}}{\sqrt{|H_n(\mathrm{i}\rho)|^2-|H_n(0)|^2}}}, & \hbox{$n$ even.}
           \end{array}
         \right.
\end{equation}
Moreover, for $n\gg1$, we have the following more explicit error bound
\begin{equation}\label{eq:errorHSI2}
\|f - h_{n} \|_{L^{\infty}(\mathbb{R})} \leq \mathcal{K} n^{1/4} \mathrm{e}^{-\rho\sqrt{2n}},
\end{equation}
where $\mathcal{K}\sim \mathrm{e}^{\rho^2/2}\widehat{V}/(2^{1/4}\pi^{3/4}\rho)$. 
\end{theorem}
\begin{proof}
From \eqref{eq:hn} we see that $\mathrm{e}^{x^2/2}h_{n}(x)$ is a polynomial of
degree $n-1$ which interpolates $\mathrm{e}^{x^2/2}f(x)$ at the points
$\{x_j\}_{j=1}^{n}$. Combining this observation with Lemma
\ref{lem:ContourInterp} gives
\[
\mathrm{e}^{x^2/2}f(x) - \mathrm{e}^{x^2/2}h_{n}(x) =
\frac{1}{2\pi\mathrm{i}} \int_{\partial{\mathcal{S}}_{\rho}}
\frac{H_{n}(x) \mathrm{e}^{z^2/2}f(z)}{H_{n}(z)(z-x)} \mathrm{d}z,
\]
or equivalently,
\[
f(x) - h_{n}(x) =
\frac{\sqrt{\gamma_{n}} }{2\pi\mathrm{i}} \psi_{n}(x) \int_{\partial{\mathcal{S}}_{\rho}}
\frac{\mathrm{e}^{z^2/2}f(z)}{H_{n}(z)(z-x)} \mathrm{d}z = \frac{1}{2\pi\mathrm{i}}  \int_{\partial{\mathcal{S}}_{\rho}}
\frac{\psi_{n}(x) f(z)}{\psi_{n}(z) (z-x)} \mathrm{d}z.
\]
Combining this with Lemma \ref{lem:MinHerm}, \eqref{eq:PsiBound} and the fact that $\min|z-x|=\rho$ for $z\in\partial{\mathcal{S}}_{\rho}$ and $x\in\mathbb{R}$, we immediately obtain \eqref{eq:errorHSI1}. As for \eqref{eq:errorHSI2}, it follows by combining \eqref{eq:errorHSI1} with \eqref{eq:HermFunAsy} and the duplication formula and ratio asymptotics for the gamma functions (see \cite[Equations (5.5.5) and (5.11.13)]{Olver2010}). This ends the proof.
\end{proof}

\section{Extensions}\label{sec:Extension}
In this section, we extend our analysis to two topics that are of interest in quadrature theory and spectral methods.

\subsection{Gauss--Hermite quadrature}
Gaussian quadrature formulas play a central role in quadrature theory and are widely used in scientific computing. Consider the following integral and its approximation
\begin{equation}
I(f) = \int_{\mathbb{R}} \mathrm{e}^{-x^2} f(x) \mathrm{d}x \approx
\sum_{k=1}^{n} w_k f(\tau_k) := Q_n(f),
\end{equation}
where $\{\tau_k\}$ and $\{w_k\}$ are the nodes and weights of the quadrature $Q_n(f)$, respectively. It is well-known that $Q_n(f)$ is the Gauss--Hermite quadrature
whenever $I(f)=Q_n(f)$ for $f\in\mathbb{P}_{2n-1}$, which achieves the maximal order of exactness for polynomials.

The convergence rate of Gauss--Hermite quadrature for analytic integrand has been studied in the past few decades. Barrett in \cite{Barrett1961} mentioned the root-exponential convergence for functions that are analytic in an infinite strip with a suitable decay at infinity. However, neither explicit proof nor the decay condition for which this rate holds were given therein\footnote{Davis and Rabinowitz in their classical monograph \cite[Equation~(4.6.1.18)]{Davis1984} also quoted Barrett's result on the root-exponential convergence of Gauss--Hermite quadrature without explicit conditions on the behavior of the integrand at infinity, but just mentioned ``certain size conditions at $z=\infty$, not mentioned here, must also be met by the integrand''.}. Donaldson and Elliott in \cite{Donaldson1972} established a contour integral representation for the remainder of Gauss--Hermite quadrature. Unfortunately, they never discussed the convergence rate of Gauss--Hermite quadrature.  Recently, Xiang in \cite{Xiang2012b} established the exponential convergence of Gauss--Hermite quadrature for entire functions whose Taylor coefficients decay supergeometrically.

In what follows, by using the contour integral representation of the remainder for Hermite spectral interpolations given in Lemma \ref{lem:ContourInterp}, we shall prove the root-exponential convergence of Gauss--Hermite quadrature for analytic integrands that subject to explicit restrictions on the behavior at infinity.
\begin{theorem}\label{thm:estGH}
If $f$ is analytic in the strip $\mathcal{S}_{\rho}$ for
some $\rho>0$ and $|f(z)|\leq\mathcal{K}|z|^{\sigma}$ for some
$\sigma\in\mathbb{R}$ as $|z|\rightarrow\infty$ within the strip and if
$V<\infty$, then for $n\geq\{\lfloor\sigma\rfloor+1,0\}$ we have
\begin{equation}\label{eq:estGH0}
|I(f) - Q_n^{\mathrm{GH}}(f)| \leq V \left\{
           \begin{array}{ll}
           {\displaystyle \left|\frac{\Phi_{n}(\mathrm{i}\rho)}{H_{n}(\mathrm{i}\rho)}\right|},          & \hbox{$n$ odd,} \\[3ex]
           {\displaystyle \frac{|\Phi_{n}(\mathrm{i}\rho)|}{\sqrt{|H_{n}(\mathrm{i}\rho)|^2 - |H_n(0)|^2}}}, & \hbox{$n$ even,}
           \end{array}
         \right.
\end{equation}
where $Q_n^{\mathrm{GH}}$ stands for the $n$-point Gauss--Hermite quadrature. Moreover, for $n\gg1$, we have the following more explicit error bound
\begin{equation}\label{eq:estGH}
|I(f) - Q_n^{\mathrm{GH}}(f)| \leq \mathcal{K} \mathrm{e}^{-2\rho\sqrt{2n}},
\end{equation}
and $\mathcal{K}\sim\mathrm{e}^{\rho^2} V$.
\end{theorem}
\begin{proof}
Let $p_{n}$ be the polynomial of degree $n-1$ determined through \eqref{eq:HermiteInter}. Recall that Gauss--Hermite quadrature is an interpolatory quadrature rule, it follows that
\begin{align}
I(f) - Q_n^{\mathrm{GH}}(f) &= \int_{\mathbb{R}} \mathrm{e}^{-x^2} (f(x) - p_{n}(x))
\mathrm{d}x. \nonumber
\end{align}
This, together with Lemma \ref{lem:ContourInterp}, implies that
\begin{align}
I(f) - Q_n^{\mathrm{GH}}(f) &= \int_{\mathbb{R}} \frac{\mathrm{e}^{-x^2}}{2\pi\mathrm{i}} \int_{\partial{\mathcal{S}}_{\rho}}
\frac{H_{n}(x)f(z)}{H_{n}(z)(z-x)} \mathrm{d}z \mathrm{d}x = \int_{\partial{\mathcal{S}}_{\rho}} \frac{\Phi_{n}(z)}{H_{n}(z)} f(z) \mathrm{d}z, \nonumber
\end{align}
where we have exchanged the order of integration in the last step. It is then readily seen that
\begin{align}
|I(f) - Q_n^{\mathrm{GH}}(f)| &\leq \max_{z\in\partial{\mathcal{S}}_{\rho}} \left|\frac{\Phi_{n}(z)}{H_{n}(z)} \right| \int_{\partial{\mathcal{S}}_{\rho}} |f(z)| |\mathrm{d}z| = V \max_{z\in\partial{\mathcal{S}}_{\rho}} \left|\frac{\Phi_{n}(z)}{H_{n}(z)} \right|. \nonumber
\end{align}
The desired result \eqref{eq:estGH0} follows by combining the last bound with Lemmas \ref{lem:MaxPhi} and \ref{lem:MinHerm}. As for \eqref{eq:estGH}, it follows by combining \eqref{eq:estGH0} with \eqref{eq:HermFunAsy}, \eqref{eq:PhiFormula1} and \eqref{eq:AsyU} and we omit the details.
\end{proof}

\begin{remark}
As far as we know, a proof on the root-exponential convergence rate of Gauss-Hermite quadrature for analytic functions is still lacking. Here we fills up this gap. Moreover, from a recent result on the worst-case error bound of Gauss-Hermite quadrature in \cite{Goda2024}, we can conclude that the root-exponential convergence rate of Gauss-Hermite quadrature is actually sharp.
\end{remark}

\begin{figure}[htbp]
\centering
\includegraphics[width=0.6\textwidth]{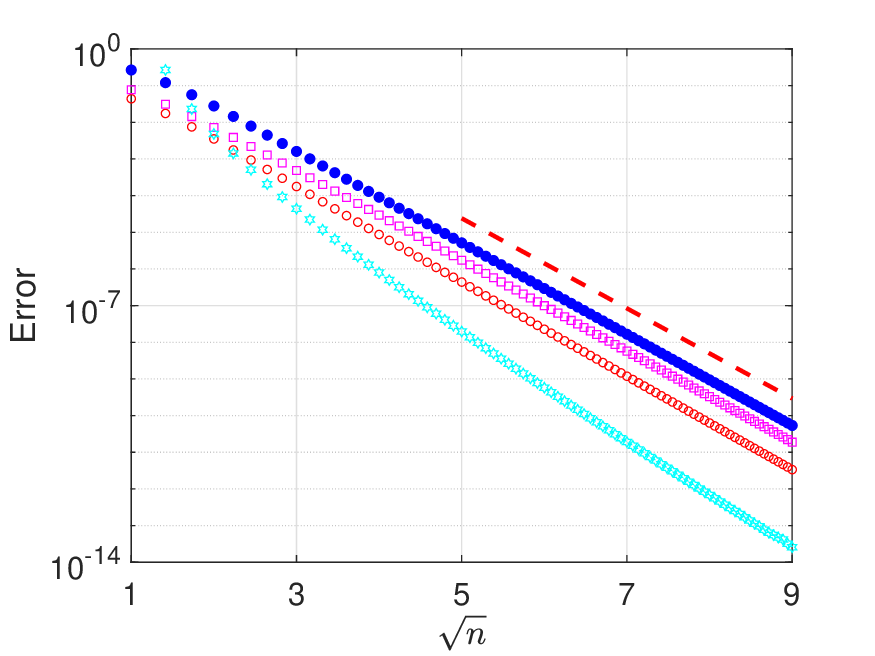}
\caption{The error of $Q_n^{\mathrm{GH}}(f)$ as a function of $\sqrt{n}$ for $f(x)=\mathrm{e}^{-x^2}/\sqrt{1+x^2}$ (dots), $f(x)=1/\sqrt{1+x^2}$ (boxes), $f(x)=\log(1+x^2)/(4+x^2)$ (circles) and $f(x)=(1+x^2)^{5/2}$ (hexagrams). The dashed line shows the rate $O(\mathrm{e}^{-2\sqrt{2n}})$.}\label{fig:GH}
\end{figure}

In Figure \ref{fig:GH} we plot the errors of Gauss--Hermite quadrature for the functions
\[
f(x) = \frac{\mathrm{e}^{-x^2}}{\sqrt{1+x^2}}, \quad \frac{\log(1+x^2)}{4+x^2}, \quad \frac{1}{\sqrt{1+x^2}}, \quad (1+x^2)^{5/2}.
\]
It is easily seen that these functions are all analytic in the strip $\mathcal{S}_{\rho}$ with $\rho=1-\epsilon$, where $\epsilon>0$ is arbitrarily close to zero, but have different asymptotic behaviors at infinity. Moreover, it is easily checked that the first two functions correspond to $V<\infty$ and the last two functions correspond to $V=\infty$, and therefore, by Theorem \ref{thm:estGH}, the associated rates of convergence of Gauss--Hermite quadrature for the first two functions are $O(\mathrm{e}^{-2\sqrt{2n}})$. We can see from Figure \ref{fig:GH} that the actual rates of Gauss--Hermite quadrature for the first two functions are consistent with our theoretical prediction. As for the last two functions, we see that Gauss--Hermite quadrature for $f(x)=1/\sqrt{1+x^2}$ still converges at the rate $O(\mathrm{e}^{-2\sqrt{2n}})$, which implies that the conditions of Theorem \ref{thm:estGH} can be relaxed further, while for $f(x)=(1+x^2)^{5/2}$ Gauss--Hermite quadrature converges at a slightly faster rate.

\begin{remark}
Gaussian quadrature rules are known to be optimal in the sense that they achieve the maximal degree of exactness. However, the non-optimality of Gauss-Hermite quadrature has been recognized by Curbera in \cite{Curbera1998} from the perspective of worst-case error of quadrature rules for Lipschitz functions and further analyzed for analytic functions in \cite{Goda2024,Trefethen2022,Trefethen2022b}. On the other hand, the evaluation of $I(f)$ by trapezoidal rule, Gauss-Legendre and Clenshaw-Curtis quadrature on a truncated interval was also extensively studied (see, e.g., \cite{Goodwin1949,Stenger1981,Sugihara1997,Trefethen2022}). For example, under the assumptions that $f(x)$ is analytic inside a strip $\mathcal{S}_{\rho}$ and satisfies certain restrictions on the asymptotic behavior at infinity within the strip, Sugihara analyzed the convergence rate of trapezoidal rule in \cite{Sugihara1997} and Trefethen analyzed the convergence rate of Gauss-Legendre and Clenshaw-Curtis quadrature in \cite{Trefethen2022}. Both of their results showed that such an $(n+1)$-point method on the interval $[-n^{1/3},n^{1/3}]$ converges at the rate $O(\exp(-Cn^{2/3}))$ for some $C>0$ as $n\rightarrow\infty$. Although the convergence rates established therein are faster than that of Gauss--Hermite quadrature for large $n$, we still found some examples, e.g., $f(x)=1/(x^2+4)$, which is not a polynomial, indicating that the latter method gives more accurate result than the former when $n$ is not large \cite{Trefethen2024}. For more discussions on sharp convergence rates and lower bounds of trapezoidal rule, Gauss-Legendre and Clenshaw-Curtis quadrature for integrals of the form $\int_{\mathbb{R}} f(x) \mathrm{d}x$, where $f(x)$ is analytic in some strip region and satisfies various decay rates at infinity, we refer to \cite{Goda2024,Sugihara1997}.
\end{remark}

\subsection{The scaling factor of Hermite approximation}
\begin{figure}[htbp]
\centering
\includegraphics[width=0.49\textwidth]{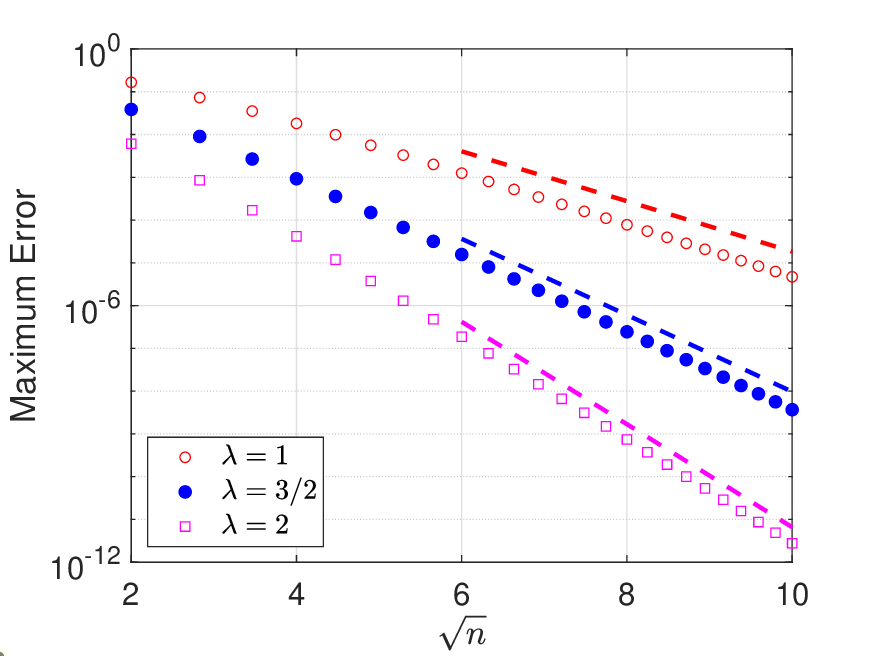}
\includegraphics[width=0.49\textwidth]{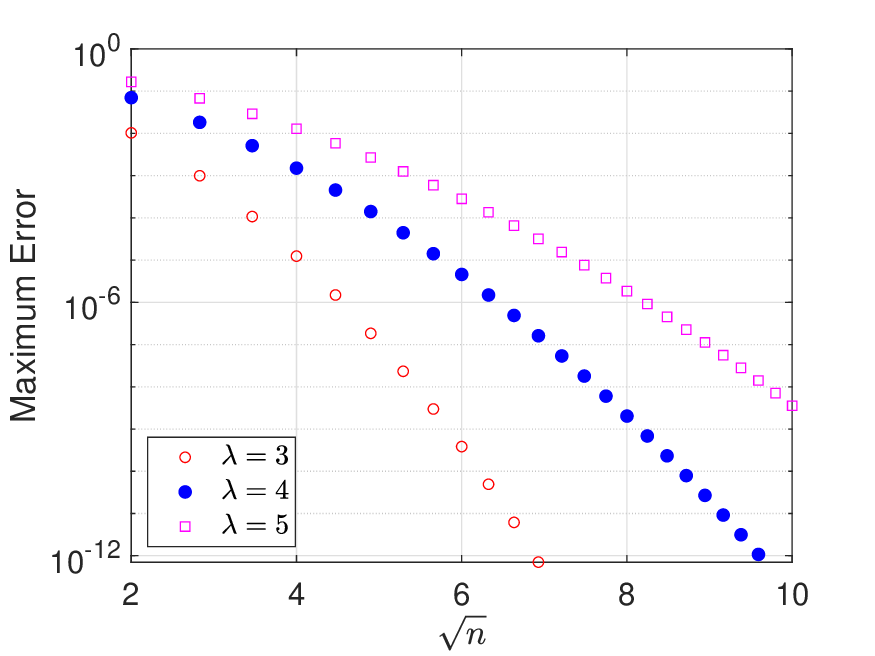}
\caption{Maximum error of $f_{n,\lambda}^{\mathrm{SF}}(x)$ for $f(x)=\mathrm{e}^{-2x^2}/(x^2+1)$ as a function of $\sqrt{n}$ for several different scaling factors. The dashed lines in the left panel show the predicted rate $O(n^{1/4}\mathrm{e}^{-\lambda\sqrt{2n}})$.}\label{fig:scaling}
\end{figure}

In practice, it is common to use the scaled Hermite functions as the basis functions in order to accelerate the rate of convergence; cf. \cite{Boyd2000,Tang1993}. 
More specifically, we consider the Hermite approximation of the form
\begin{equation}
f_{n,\lambda}^{\mathrm{SF}}(x) = \sum_{k=0}^{n} \alpha_k \psi_k(\lambda x), \qquad  \alpha_k
= \lambda\int_{\mathbb{R}}f(x)\psi_k(\lambda x) \mathrm{d}x,
\end{equation}
where $\psi_k$ is the Hermite function defined in \eqref{def:HermFunc} and $\lambda>0$ is a scaling factor. It comes out that our convergence analysis is also useful to some extent in finding the scaling factor. If $f(z)$ is analytic inside and on the strip $\mathcal{S}_{\rho}$ for some $\rho>0$ and
$$
\int_{\partial\mathcal{S}_{\rho}}|\mathrm{e}^{(\lambda z)^2/2}f(z)| |\mathrm{d}z| < \infty,
$$
we then obtain from Theorem \ref{thm:ProHermFunc} that $f_{n,\lambda}^{\mathrm{SF}}(x)$ converges at the rate $O(n^{1/4}\mathrm{e}^{-\lambda\rho\sqrt{2n}})$ in the maximum norm. Clearly, the scaling factor $\lambda$ in this case should be chosen as large as possible to maximize the convergence rate. As an example, we consider the function $f(x)=\mathrm{e}^{-2x^2}/(x^2+1)$. Our analysis implies that the convergence rate of $f_{n,\lambda}^{\mathrm{SF}}(x)$ in the maximum norm is $O(n^{1/4}\mathrm{e}^{-\lambda\sqrt{2n}})$ whenever $\lambda\leq2$. In the left panel of Figure \ref{fig:scaling} we plot the maximum errors of $f_{n,\lambda}^{\mathrm{SF}}(x)$ for $\lambda=1,3/2,2$. As expected, the convergence rate of $f_{n,\lambda}^{\mathrm{SF}}(x)$ is indeed $O(n^{1/4}\mathrm{e}^{-\lambda\sqrt{2n}})$ and $\lambda=2$ achieves a much faster convergence rate than $\lambda=1$. When $n=100$, we see that the accuracy of $f_{n,\lambda}^{\mathrm{SF}}(x)$ has been improved dramatically from $10^{-6}$ for $\lambda=1$ to $10^{-12}$ for $\lambda=2$. Nevertheless, we point out that $\lambda=2$ is still not the optimal scaling factor. As illustrated in the right panel of Figure \ref{fig:scaling}, it is easily seen that the choice $\lambda=3$ achieves a faster convergence rate than $\lambda=2$, and the convergence rate is actually faster than root-exponential.

\section{Conclusion}\label{sec:Conclusion}
In this work we have presented a rigorous convergence analysis of Hermite spectral approximations for functions that are analytic within an infinite strip and satisfies certain restrictions on the asymptotic behavior at infinity within the strip. By exploring some remarkable contour integral representations for the Hermite coefficients and the remainder of Hermite spectral interpolations, we derive some sharp error bounds for Hermite approximations in the weighted and maximum norms. Extensions of our analysis to Gauss--Hermite quadrature and the scaling factor of Hermite approximation are also discussed. In particular, we derive the root-exponential convergence of Gauss--Hermite quadrature for analytic integrands and provide explicit conditions under which this rate holds true.

Finally, it is worthwhile to point out that some issues related to our work still remain. For example, the conditions of Theorems \ref{thm:ProHermFunc} and \ref{thm:estGH} might be further relaxed. Moreover, Figure \ref{fig:scaling} implies that the choice of the optimal scaling factor still does not follow from the current analysis. We leave theses issues for future studies.

\section*{Acknowledgements}
The authors wish to thank the anonymous referees for their valuable comments to improve the presentation of the paper. Haiyong Wang would like to thank Simon Chandler-Wilde, David Hewett, Daan Huybrechs, Arieh Iserles and Nick Trefethen for helpful discussions and feedbacks of the present work during the workshop ``Singular and Oscillatory Integration: Advances and Applications" held at University College London from June 24 to June 26, 2024 and Weizhu Bao, Zhongqing Wang and Zhimin Zhang for advice and suggestions. His research was supported in part by the National Natural Science Foundation of China under grant number 12371367 and by the Hubei Provincial Natural Science Foundation of China under grant number 2023AFA083. Lun Zhang was supported in part by the National Natural Science Foundation of China under grant number 11822104 and by ``Shuguang Program'' supported by Shanghai Education Development Foundation and Shanghai Municipal Education Commission.

\end{document}